\documentclass[11pt, reqno] {amsart}

\usepackage{graphicx,amsthm}
\usepackage{verbatim,amsmath,amscd,amssymb,color}
\usepackage{cases}
\usepackage{empheq}
\usepackage[mathscr]{eucal}
\usepackage{xcolor}
\usepackage{verbatim}
\usepackage{tikz}
\usepackage{mathrsfs} 

\colorlet{mycol}{black}
\usetikzlibrary{positioning,arrows.meta,backgrounds}
\usepackage{txfonts}
\definecolor{navyblue}{RGB}{20,20,140}

\usepackage{epic,eepic,multicol}
\numberwithin{equation}{section}
\input xy
\xyoption{all}

\def\m@th{\mathsurround=0pt}
\def\fsquare(#1,#2){
\hbox{\vrule$\hskip-0.4pt\vcenter to #1{\normalbaselines\m@th
\hrule\vfil\hbox to #1{\hfill$\scriptstyle #2$\hfill}\vfil\hrule}$\hskip-0.4pt
\vrule}}

\newcommand{\cA}{{\mathcal A}}
\newcommand{\cB}{{\mathcal B}}

\newcommand{\cK}{{\mathcal K}}

\newcommand{\cT}{{\mathcal T}}

\newcommand{\frg}{\mathfrak g}
\newcommand{\frh}{\mathfrak h}

\newcommand{\frn}{\mathfrak n}

\newcommand{\bbA}{\mathbb A}
\newcommand{\bbB}{\mathbb B}

\newcommand{\bbQ}{\mathbb Q}

\newcommand{\bbZ}{\mathbb Z}

\newcommand{\bfE}{\mathbf E}

\newcommand{\scC}{{\mathscr C}}
\newcommand{\scD}{{\mathscr D}}

\newcommand{\ch}{\mathrm{ch}}

\newcommand{\op}{\mathrm{op}}

\newcommand{\wt}{\mathrm{wt}}

\newcommand{\lbr}{\begin{bmatrix}}
\newcommand{\rbr}{\end{bmatrix}}
\newcommand{\cd}{commutative diagram }

\def\ge{\frg}

\def\al{\alpha}
\def\aq{\cA_q(\frn)}

\def\bfk{{\bf k}}
\def\ci(#1,#2){c_{#1}^{(#2)}}
\def\Ci(#1,#2){C_{#1}^{(#2)}}
\def\mpp(#1,#2,#3){#1^{(#2)}_{#3}}
\def\bCi(#1,#2){\ovl C_{#1}^{(#2)}}
\def\ch(#1,#2){c_{#2,#1}^{-h_{#1}}}
\def\cc(#1,#2){c_{#2,#1}}

\def\bfi{{\mathbf i}}

\def\beneme{\begin{enumerate}}
\def\beq{\begin{equation}}
\def\beqn{\begin{eqnarray}}
\def\beqnn{\begin{eqnarray*}}

\def\bbra#1,#2,#3{\left\{\begin{array}{c}\hspace{-5pt}
#1;#2\\ \hspace{-5pt}#3\end{array}\hspace{-5pt}\right\}}
\def\cd{\cdots}
\def\ci(#1,#2){c_{#1}^{(#2)}}

\def\del{\delta}
\def\Del{\Delta}

\def\eit{\tilde{e}_i}
\def\Eit{\widetilde{E}_i}
\def\eneme{\end{enumerate}}

\def\eeq{\end{equation}}
\def\eeqn{\end{eqnarray}}
\def\eeqnn{\end{eqnarray*}}
\def\fit{\tilde{f}_i}
\def\Fit{\widetilde{F}_i}

\def\gau#1,#2{\left[\begin{array}{c}\hspace{-5pt}#1\\
\hspace{-5pt}#2\end{array}\hspace{-5pt}\right]}

\def\ify{\infty}

\def\lan{\langle}

\def\max{{\rm max}}
\def\lm{\lambda}
\def\Lm{\Lambda}
\def\mapright#1{\smash{\mathop{\longrightarrow}\limits^{#1}}}

\def\nd{\noindent}
\def\nn{\nonumber}

\def\ot{\otimes}
\def\op{\oplus}

\def\ovl{\overline}

\def\qq{\qquad}
\def\q{\quad}
\def\qed{\hfill\framebox[2mm]{}}

\def\ran{\rangle}
\def\rgmod{\hbox{$R$-\hbox{gmod}}}

\def\til{\tilde}
\def\tilRgmod{\widetilde{R}\hbox{-gmod}}
\def\tilRgmod{\hbox{$\widetilde{R\hbox{-gmod}}$}}

\def\uq{U_q(\ge)}

\def\uqm{U^-_q(\ge)}

\def\TY(#1,#2,#3){#1^{(#2)}_{#3}}

\def\vep{\varepsilon}
\def\vp{\varphi}

\def\xxi(#1,#2,#3){\displaystyle {}^{#1}\Xi^{(#2)}_{#3}}

\def\wtil{\widetilde}

\def\m@th{\mathsurround=0pt}
\def\fsquare(#1,#2){
\hbox{\vrule$\hskip-0.4pt\vcenter to #1{\normalbaselines\m@th
\hrule\vfil\hbox to #1{\hfill$\scriptstyle #2$\hfill}\vfil\hrule}$\hskip-0.4pt
\vrule}}

\theoremstyle{definition}
\newtheorem{df}{Definition}[section]
\newtheorem{thm}[df]{Theorem}
\newtheorem*{introthm}{Theorem 7.4}
\newtheorem{pro}[df]{Proposition}
\newtheorem{lem}[df]{Lemma}
\newtheorem{ex}[df]{Example}
\newtheorem{cor}[df]{Corollary}

\title[Characterization of the unit object in localized quantum unipotent category]
{ Characterization of the unit object in localized quantum unipotent category }

\author{ K\textsc{oh} M\textsc{atsuura}}\thanks{K.M.;
{Division of Mathematics, 
Sophia University, Kioicho 7-1, Chiyoda-ku, Tokyo 102-8554,
Japan}\\
Email:\texttt{k-matsuura-4w0@eagle.sophia.ac.jp}\,\,}
\author{ T\textsc{oshiki} N\textsc{akashima}}\thanks{T.N.;
{Division of Mathematics, 
Sophia University, Kioicho 7-1, Chiyoda-ku, Tokyo 102-8554,
Japan}\\
Email:\texttt{toshiki@sophia.ac.jp},\,\,
T.N is supported in part by
JSPS Grants in Aid for Scientific Research $\#$25K06939. }

\date{}

\begin{document}
\maketitle
\begin{abstract}
For the quiver Hecke algebra $R$, let $\rgmod$ be the category of finite-dimensional 
graded $R$-modules, and let $\wtil\rgmod[w]$ be the localization of  $\rgmod$. Kashiwara and the second author showed the set of equivalence classes of simple objects up to grading shifts $\mathrm{Irr}(\wtil\rgmod[w])$ in $\wtil\rgmod[w]$ has a crystal structure, and $\mathrm{Irr}(\wtil\rgmod[w])$ is isomorphic to the so-called cellular crystal $\bbB_\bfi$. 
This isomorphism induces a function $\vep_i^*$ on $\bbB_\bfi$. We give an explicit formula of $\vep_i^*$, and using this formula, we give a characterization of the unit object of $\wtil\rgmod[w_0]$ for the case of classical finite types.

\end{abstract}

\tableofcontents

\section{Introduction}


The quiver Hecke algebra has been introduced independently by 
Khovanov-Lauda (\cite{K-L,K-L2}) and Rouquier (\cite{Rou}). It is a family of $\bbZ$-graded ${\bf k}$-algebras
$R(\beta)$ indexed by $\beta=\sum_im_i\al_i\in Q_+$, where $Q_+$ is the 
positive root lattice and ${\bf k}$ is a field.
For a symmetrizable Kac-Moody Lie algebra $\ge$, 
one of the remarkable properties of the algebra $R:=\bigoplus_{\beta\in Q_+}R(\beta)$ is that 
it categorifies the nilpotent half of the quantum algebra $\uqm$ and 
the unipotent quantum coordinate ring $\aq$. Specifically,
there exist isomorphisms of algebras:
\[
\cK(R\hbox{\rm-proj})\cong U^-_q(\ge)_{\bbZ[q,q^{-1}]},\qq\qq\cK(\rgmod)\cong 
\cA(\mathfrak n)_{\bbZ[q,q^{-1}]},
\]
where $\cK(R\hbox{\rm-proj})$ (resp. $\cK(\rgmod)$) is the Grothendieck ring of 
$R\hbox{\rm-proj}$ (resp. $\rgmod$) the category of graded projective (resp. finite-dimensional) $R$-modules equipped with the convolution product as multiplication.

 Crystals are combinatorial objects developed from the theory of the crystal bases of the quantum group \cite{K3}. A crystal consists of a 6-tuple $(B,\wt, \{\vep_i\},\{\vp_i\}, \{\eit\},\{\fit\})_{i\in I}$
with a set $B$, maps $\eit,\fit:B\sqcup\{0\}\to B\sqcup\{0\}$,
and functions $\vep_i,\vp_i:B\to \bbZ\sqcup\{-\ify\}$, $\wt:B\to P$ satisfying certain conditions (see Definition \ref{def of crystal}). Crystals admit a natural tensor product structure \cite{K3}.
An important example is given by the tensor product of crystals $B_i$ associated with a reduced word $i_1i_2\cdots i_l$ where $B_i$ is in Example \ref{e.g.cr2} below.  Owing to the existence of special isomorphisms $\phi_{ij}^{(k)}$  ($k=0,1,2,3$), called braid-type isomorphisms, $B_{i_1}\ot B_{i_2}\ot\cdots \ot B_{i_l}$ does not depend on the choice of a reduced word $\bfi=i_1i_2\cdots i_l$ i.e., for two reduced words, $i_1\cdots i_l$ and $j_i\cdots j_l$, we have an isomorphism $B_{i_1}\ot\cdots\ot B_{i_l}\cong B_{j_1}\ot\cdots\ot B_{j_l}.$ Thus, let us  $\bbB_\bfi$ denote this crystal and call it a cellular crystal associated to $\bfi$ (see Proposition \ref{pro-braid} and Definition \ref{def of cellular}). Another fundamental example of crystal is the crystal basis $B(\ify)$ of the nilpotent half of the quantum algebra $\uqm$. Since, as we mentioned,  the quiver Hecke algebra categorifies the nilpotent half of the quantum algebra,  it is natural to expect that it also categorifies its crystal structure. Indeed, Lauda and Vazirani \cite{L-V} showed that the family $\bbB$ of finite dimensional  self-dual simple $R$-modules has a crystal structure and is isomorphic to the crystal $B(\ify)$, which can be viewed as a categorical realization of the crystal $B(\ify)$.

In a series of articles \cite{Loc, Loc2, Loc3}, the localization of  monoidal categories of $R$-modules is studied.
For a graded monoidal category $(\cT=\oplus_{\lm\in \Lm}\cT_\lm,\ot)$, where 
$\Lm$ is a $\bbZ$-lattice, a real commuting family of graded braiders $\{(C_i,R_{C_i},\phi_i)\}_{i\in I}$ consists of objects $C_i$ in $\cT$, morphisms $R_{C_i}:(C_i\ot -)\to (-\ot C_i)$ and $\bbZ$-valued functions $\phi_i$ on $\Lm$ satisfying certain conditions. Then, we obtain a localization $\wtil\cT$ of the category $\cT$ which enjoys nice properties as in Proposition \ref{til-pro} and a canonical functor $\Phi: \cT \to \wtil\cT$ such that the objects $\Phi(C_i)$ are invertible, the morphism $\Phi(R_{C_i}): \Phi(C_i)\ot\Phi(X) \to \Phi(X)\ot\Phi(C_i)$ are isomorphisms for all $i$ and all $X \in\cT$. Here, the real commuting family of graded braiders
$\{(C_i,R_{C_i},\phi_i)\}_{i\in I}$ plays a similar role as a multiplicative set for the localization of commutative ring theory.

In $\rgmod$, there exist graded braiders $(M(w\Lambda_i, \Lambda_i), R_{M(w\Lambda_i, \Lambda_i)}, \phi_M(w\Lambda_i, \Lambda_i))$ for any $i \in I$ consisting of the determinantial module $M(w\Lambda_i, \Lambda_i)$ and the morphism induced from $R$-matrices
\[
R_{M(w\Lambda_i, \Lambda_i)}: M(w\Lambda_i, \Lambda_i)\circ X \to X\circ M(w\Lambda_i, \Lambda_i) 
\] 
for any $X\in \rgmod$ of homogeneous degree $-(w\Lambda_i+\Lambda_i, \mathrm{wt}(X))$. Thus, we obtain the localizaion $\wtil\rgmod[w]:=\rgmod[M(w\Lambda_i, \Lambda_i)^{\circ-1} ; i \in I]$ of the category $\rgmod$.
Applying this localization method to the subcategory $\scC_w$ of $R$-gmod associated with a Weyl group element $w\in W$, one has the localized category $\wtil\scC_w$. In particular, for finite type $\ge$ and the longest element $w=w_0$  in $W$, $\scC_{w_0}$ coincides with the category $R$-gmod. Hence, we obtain $\wtil\rgmod[w_0] = \wtil\scC_{w_0}$. 

Since $\rgmod$ possesses the crystal structure, it is natural to expect $\wtil\scC_w$ also possess a crystal sturcture.
In \cite{KN}, it is shown that the set of equivalence clasess of simple objects in  $\mathrm{Irr}(\wtil\scC_w)$ possesses a crystal structure and is isomorphic to the cellular crystal $\bbB_\bfi:=B_{i_1}\ot\cd\ot B_{i_k}$ where $w=s_{i_1}s_{i_2}\cdots s_{i_k}$ is a reduced expression of $w$.  The isomorphism is given explicitly by
\begin{align}
\text{Irr}(\wtil\scC_w) \longrightarrow \bbB_\bfi \q \text{by}  \q \wtil C_{\lambda}^{-1}\circ Q_w(M) \longmapsto K_{\bfi}(M)-K_{\bfi}(C_{\lambda})
\end{align}
where $K_{\bfi}(M)=(c_1, \dots c_l)$, $c_k=\vep_{i_k}^*(M_k)$, $M_k= (\wtil E^*_{i_k})^{c_k}(M_k)$ for $1 \leq k \leq l$ and $M_l=M.$ Furthermore, for $w \in W$ and $i \in I$ satisfiying $ w'=ws_i < w$, this isomorphism is obtained inductively by
\begin{align}\label{inductiveisom}
\mathrm{Irr}(\wtil\scC_w)\to \mathrm{Irr}(\wtil\scC_{w'})\times \bbZ \quad (X \mapsto (\bfE_{w',w}(X), \vep_i^*(X))
\end{align}
(see Section 4 for more details). As mentioned above, for finite type $\ge$ and the longest element $w=w_0$  in $W$,
we have $\mathrm{Irr}(\wtil\rgmod[w_0])\cong \bbB_\bfi$. Therefore, this result can be regarded as the generalization of the result in \cite{N}.

In this article, we restrict ourselves to the case where $\ge$ is classical finite type  and the longest element $w_0$. Since $w_0 \in W$ is the longest element, any index $i \in I$ can be moved to the rightmost by applying braid moves to $\bfi$ (see Definition \eqref{def of braid-moves} for the definition of the braid moves). This procedure induces an isomorphism of the crystals
\begin{align}
\phi : \bbB_\bfi \to \bbB_{\mathbf{i'}}\ot B_i \quad (x \mapsto \phi(x)=b'\ot (a)_i)
\end{align}
where $\bfi$ and $\mathbf{i'}i$ are the reduced words of $w_0$ and $\phi$ is the composition of the braid type isomorphisms (see Definition \ref{pro-braid} for the definition of braid type isomorphisms).

Let us present the main results of this article. We give a explicit description of $\vep_i^*$ for a fixed longest word $\bfi$ of classical finite types. For example, for type $A_n$, let $\bfi=(1)(21) \dots (n-1 \dots 21)(n\dots21)$ be a reduced word, and set
\[
x=(z_{1,1})_1\ot(z_{1,2})_2\ot(z_{2,1})_1\ot\cdots\ot(z_{1,n})_n\ot\cdots\ot(z_{n-1,2})_2\ot(z_{n,1})_1 \in\bbB_{\bfi}
\]
where the double indices $(j, i_k)$ of $z_{j, i_k}$ indicates $j$-th occurrence of $i_k$. Then
$\vep_{i}^*(x)$ is given by
\begin{align*}
&\vep_{i}^*(x)=\underset{1 \leq k \leq i}{\max}(z_{n-i+2, k-1}- z_{n-i+1, k}) \q (1\leq i \leq n).
\end{align*}
where we set \ $z_{j,k}=0$, if \ $k=0$. We also obtain the explict form of $\vep_i^*(x)$ for other types in Proposition \ref{pro-vep*}.
We also provide a criterion for $X=\mathbf{1}$ in $\wtil\rgmod[w_0]$. More precisely, we obtain the following result:
\begin{introthm}
Let W be the Weyl group of classical finite type. $\bfi$ = $i_1\cd i_l$  be a reduced longest word, and $\bbB_\bfi$ be the cellular crystal associated $\textbf{i}$. For $x \in \bbB_\bfi$, the following are equivalent.
\begin{enumerate}
\item
$\mathrm{wt}(x)=0$ and for all $i \in I, \vep_{i}^*(x)=0.$
\item
$x = (0)_{i_1}\otimes \cd \otimes(0)_{i_l}\in \bbB_\bfi.$
\end{enumerate}
\end{introthm}
Since $K_{\bfi}$ is a bijection, the condition $(1)$ gives the characterization of the unit object in $\wtil\rgmod[w]$.

Although it is known that for $b\in B(\ify)$, the condition $\vep_i^*(b)=0$ for any $i \in I$ implies $b=u_{\ify}:= 1 \mathrm{mod} \ qL(\ify)$, it does not hold in the localized case. For example, for $\ge=A_3$, fix $\bfi=1(21)(321)$ as a reduced word. In this case, $\vep_i^*(x)$ is given by
\[
\vep_1^*(x)=-z_{3,1},\ \vep_2^*(x)=\max(-z_{2,1}, z_{3,1}-z_{2,2}),\ \vep_3^*(x)=\max(-z_{1,1}, z_{2,1}-z_{1,2}, z_{2,2}-z_{1,3}).
\]
Then we have $x=(1)_1\ot(1)_2\ot(0)_1\ot(1)_3\ot(1)_2\ot(0)_1$ satisfies $\vep_i^*(x)=0$ for any $i \in I$, while $x\ne u_{\ify}$.

This article is organized as follows.
In Section 2, we recall the notion of crystals and introduce the cellular crystal.
Section 3 is devoted to review the definition and basic properties of quiver Hecke algebras.
In Section 4, we briefly review the localization of a monoidal category and discuss some of its properties. We also recall the crystal structure on the localized category and the theorems establishing its isomorphism with the cellular crystal.
Finally, in Section 5, we present explicit formulae for $\vep_i^*$ and prove the main results.


\section{Crystals}
In this section, we recall the notion of crystals. Then we will introduce the cellular crystal which plays a fundamental role in this paper.
\subsection{Preliminaries}
Let $I$ be the index set. An integer-valued square matrix \ $A=(a_{ij})$ is the symmetrizable generalized Cartan matrix if it satisfies 
\begin{enumerate}
\item$a_{ii}=2$ ($i \in I$),
\item$a_{ij}\leq 0$ \ if \ $i\ne j$,
\item$a_{ij}=0$ \ if and only if \ $a_{ji}=0.$
\item there exists a diagonal matrix \ D=diag($d_i \mid i \in I$) \ with all $d_{i}$ positive integers such that DA is symmetric.
\end{enumerate}

The dual weight lattice $P^\vee$ is a free abelian group with a basis $\{h_i \mid i\in I\}$ and let $\frh=\bbQ$ $\ot_{\bbZ}$ $P^\vee$.Then we define the weight lattice to be $P= \{ \lambda \in \frh^* \mid \lambda(P^\vee) \subset \bbZ \}$. Let $\Pi^\vee= \{h_i \mid i\in I\}$ and choose a linearly independent subset $\Pi=\{\alpha_i \mid i \in I \}$ that satisfies \ $\langle h_i, \alpha_j \rangle =\alpha_j(h_i)=a_{ij}$.We call the elements of $\Pi$ simple roots and the elements of $\Pi^\vee$ simple coroots.

\begin{df}
 A Cartan datum associated with the generalized Cartan matrix A is a quintuple $(A, P, P^\vee, \Pi, \Pi^\vee).$
\end{df}

We call $Q=\op_{i \in I}\bbZ \alpha_i$ the root lattice and $Q_{+}=\op_{i\in I} \bbZ_{\geq 0} \alpha_i$ the positive root lattice. For an element  $\beta=\sum_i m_i\al_i\in Q_+$, we define $|\beta|=\sum_i m_i$, which is called the height of $\beta$. We also define a symmetric bilinear form ( \ ,\ ) on $\frh^*$ satisfying $(\alpha_i,\alpha_j)=d_ia_{ij}$ \ $(i,j \in I)$ and $\langle h_i, \lambda \rangle =$ $2(\alpha_i, \lambda)\over{(\alpha_i, \alpha_i)}$ for any $\lambda \in \frh^*$ and  $i\in I.$

Let $\uq$ be the quantum group associated with $(A, P, P^\vee, \Pi, \Pi^\vee)$, which is a $\bbQ$-algebra generated by $e_i$, $f_i \ (i \in I)$ and $q^h \ ( h\in P^\vee)$ with certain relations. Let $\uqm$ denote the subalgebra of $\uq$ generated by $f_i \ (i \in I)$.





Let $W=\lan s_i \ran_{i\in I}$ be the Weyl group associated with $P$,
where $s_i$ is the simple reflection defined by $s_i(\lm)
=\lm-\lan h_i,\lm\ran\al_i$ $(\lm\in P)$. In the sequel, we assume $A=(a_{ij})$ and $W$ are of finite type.
For the numbering of Dynkin diagrams of Lie algebras, we adopt the ones in \cite[Chap.4]{Kac}.
\subsection{Language of crystal}
\begin{df}\label{def of crystal}
A 6-tuple $(B,\wt, \{\vep_i\},\{\vp_i\}, \{\eit\},\{\fit\})_{i\in I}$ is 
a {\it crystal} if $B$ is a set and there exists an element denoted by $0\not\in B$ and maps:
\begin{eqnarray*}
&&{\rm wt}:B\to P,\label{wtp-c}\q
\vep_i:B\to \bbZ\sqcup\{-\ify\},\q\vp_i:B\to \bbZ\sqcup\{-\ify\}\q\,(i\in I)
\\
&&\eit:B\sqcup\{0\}\to B\sqcup\{0\},\q
\fit:B\sqcup\{0\}\to B\sqcup\{0\}\,\,(i\in I),\label{eitfit-c}
\end{eqnarray*}
satisfying :
\begin{enumerate}
\item
$\vp_i(b)=\vep_i(b)+\lan h_i,\wt(b)\ran$.
\item
If  $b,\eit b\in B$, then $\wt(\eit b)=\wt(b)+\al_i$, 
$\vep_i(\eit b)=\vep_i(b)-1$, $\vp_i(\eit b)=\vp_i(b)+1$.
\item
If $b,\fit b\in B$, then $\wt(\fit b)=\wt(b)-\al_i$, 
$\vep_i(\fit b)=\vep_i(b)+1$, $\vp_i(\fit b)=\vp_i(b)-1$.
\item
 For $b,b'\in B$ and $i\in I$, one has $\fit b=b'$ if only if $b=\eit b'.$
\item
If $\vp_i(b)=-\ify$ for $b\in B$, then $\eit b=\fit b=0$
and $\eit(0)=\fit(0)=0$.
\end{enumerate}
\end{df}

\begin{df}
For a crystal $B$, its crystal graph is an oriented $I$-colored graph defined by
\[
\centerline{$ b\mapright{i} b'\Leftrightarrow \fit b=b'$}.
\]
\end{df}

\begin{ex}\label{ex-cry}
For $\lambda \in P$, let $T_{\lambda}=\{t_\lambda\}$ and for any  $ i \in I$, define
\[
wt(t_\lambda)=\lambda, \q \eit t_\lambda=\fit t_\lambda=0, \q \varepsilon_i(t_\lambda)=\varphi_i(t_\lambda)=-\ify.
\]Then $T_{\lambda}$ is a crystal.
\end{ex}

\begin{ex}[\cite{K1}]
There exists the crystal basis $(L(\ify),B(\ify))$ of $\uqm$ defined by
\begin{gather*}
 L(\ify):=\sum_{k\geq0,i_1,\cd,i_k\in I}\bbA\til f_{i_1}\cd\til
 f_{i_k}u_\ify,\\
B(\ify)=\{\til f_{i_1}\cd\til
 f_{i_k}u_\ify\,{\rm mod}\, q L(\ify)\,|\,k\geq0, i_1,\cd,i_k\in I\}\setminus\{0\},\\
\vep_i(b)={\rm max}\{k:\eit^kb\ne0\},\q
\vp_i(b)=\vep_i(b)+\lan h_i,\wt(b)\ran,
\end{gather*}
where $u_\ify=1\in \uq$, $\eit$ and $\fit\in{\rm End}_{\bbQ(q)}(\uqm)$ are the Kashiwara operators and $\bbA$ is a local subring of $\bbQ(q)$ at $q=0$.
\end{ex}

\begin{ex}\label{e.g.cr2}
For $i\in I$, set
$B_i:=\{(n)_i\,|\,n\in \bbZ\}$ and its crystal structure is given by 
\begin{eqnarray*}
&&\wt((n)_i)=n\al_i,\,\,
\vep_i((n)_i)=-n,\,\,
\vp_i((n)_i)=n,\,\,\label{bidata}\\
&&
\vep_j((n)_i)=\vp_j((n)_i)=-\ify\,\,\,(i\ne j),\label{biji}\\
&&
{ \eit((n)_i)=(n+1)_i,\q\fit((n)_i)=(n-1)_i},\label{bieitfit}\\
&&
\til e_j((n)_i)=\til f_j((n)_i)=0\,\,\,(i\ne j).
\label{biejfj}
\end{eqnarray*}
\end{ex}

\begin{df}
For crystals $\cB_1, \cB_2$, a {\it crystal morphism} $\Psi: \cB_1 \rightarrow \cB_2$ is a map $\Psi: \cB_1 \cup \{0\} \rightarrow \cB_2 \cup \{0\}$ such that
\begin{flushleft}
\begin{enumerate}
\item
$\Psi(0)=0$,
\item
if $b \in \cB_1 $ and $\Psi(b) \in \cB_2,$ then $wt(\Psi(b))=wt(b),$
$\vep_i(\Psi(b))=\vep_i(b)$ and $\vp_i(\Psi(b))=\vp_i(b)$ for all $i \in I$,
\item
f $b,b' \in \cB_1$, $\Psi(b), \Psi(b')\in\cB_2$ and $\fit b=b'$, then $\fit\Psi(b)=\Psi(b')$ and $\Psi(b)= \eit\Psi(b')$ for all $i\in I$.
\end{enumerate}
\end{flushleft}
$\Psi: \cB_1 \rightarrow \cB_2$ is called an embedding if $\Psi:\cB_1 \cup \{0\} \rightarrow \cB_2 \cup \{0\}$ is an injection \\and called an isomorphism if a bijection.
\end{df}

Crystals also have the tensor product, whose crystal structure is obtained by the tensor product rule of crystal bases as follows.
\begin{df}
The tensor product $\mathcal{B}_1 \otimes \mathcal{B}_2$ of crystals $\mathcal{B}_1$ and $\mathcal{B}_2$ is defined to be the set $\mathcal{B}_1 \times \mathcal{B}_2$ whose crystal structure is given by
\begin{enumerate}
\item
$wt(b_1 \otimes b_2)=\mathrm{wt}(b_1)+\mathrm{wt}(b_2),$
\item
$\varepsilon_i(b_1 \otimes b_2)=\max\{\vep_i(b_1),\vep_i(b_2)-\langle h_i, \mathrm{wt}(b_1)\rangle \},$
\item
$\vp_i(b_1\otimes b_2)=\max\{\vp_i(b_2),\vp_i(b_1)+\langle h_i, \mathrm{wt}(b_2)\rangle \},$
\item
$\eit(b_1\otimes b_2)=
\begin{cases}
\eit b_1 \otimes b_2\qq \mathrm{if}\ \vp_i(b_1)\geq \vep_i(b_2),\\b_1 \otimes \eit b_2\qq \mathrm{if}\ \vp_i(b_1)< \vep_i(b_2),
\end{cases}$
\item
$\fit(b_1 \otimes b_2)=
\begin{cases}
\fit b_1 \otimes b_2 \qq \mathrm{if}\ \vp(b_1)>\vep_i(b_2),\\b_1 \otimes \fit b_2 \qq \mathrm{if}\ \vp_i(b_1)\leq \vep_i(b_2),
\end{cases}$

\end{enumerate}
where we write $b_1\otimes b_2$ for $(b_1,b_2)\in \mathcal{B}_1 \times \mathcal{B}_2$, and we set $b_1\otimes0=0\otimes b_2=0.$

\end{df}

\subsection{Cellular crystal}
For the tensor product of the crystals $B_i$ (see Example \ref{e.g.cr2}), there exist isomorphisms, called braid type isomorphisms.
\begin{pro}[\cite{N1}]\label{pro-braid}
For $i\in I$, set $B_i=\{(z)_i\,|\,z\in \bbZ\}$ and let $A=(a_{ij})_{i,j \in I} $ be a Cartan matrix. 
Then, there exist the following isomorphisms of crystals
$\phi_{ij}^{(k)}$ 
($k=0,1,2,3$): 
\begin{align*}
&\phi_{ij}^{(0)}:B_i\ot B_j\overset{\sim}{\longrightarrow}B_j\ot B_i, \ (z_1)_i\otimes(z_2)_j \longmapsto (z_2)_j\otimes (z_1)_i, \q \text{if} \q a_{ij}a_{ji}=0, \\
&\phi_{ij}^{(1)}:B_i\ot B_j\ot B_i\mapright{\sim}B_j\ot B_i\ot B_j, \q \text{if} \q
a_{ij}a_{ji}=1, \\
&(z_1)_i\otimes(z_2)_j\otimes(z_3)_j\longmapsto (\max(z_3, z_2-z_1))_j\otimes(z_1+z_3)_i\otimes(-\max(-z_1, z_3-z_2))_j,&\\
&\phi_{ij}^{(2)}:B_i\ot B_j\ot B_i\ot B_j
\mapright{\sim}B_j\ot B_i\ot B_j\ot B_i, \q \text{if} \q
a_{ij}=-1, a_{ji}=-2, \\
&(z_1)_i \ot (z_2)_j\ot (z_3)_i \ot (z_4)_j\\
&\longmapsto (\max(z_4, z_2-2z_1, 2z_3-z_2))_j \ot (\max(z_1+z_4,z_2,z_1-z_2+2z_3))_i\\
&\ot (-\max(-z_3,-z_4-2z_1,-2z_2+2z_3-z_4))_j
\ot(-\max(-z_3+z_4, -z_1,z_3-z_2))_i.\\
&\phi_{ji}^{(2)}:B_j\ot B_i\ot B_j\ot B_i
\mapright{\sim}B_i\ot B_j\ot B_i\ot B_j, \q \text{if} \q
a_{ij}=-1, a_{ji}=-2 , \\
&(z_1)_j \ot (z_2)_i\ot (z_3)_j \ot (z_4)_i\\
&\longmapsto (\max(-z_2+z_3, -z_1+z_2, z_4))_i \ot (\max(z_1-2z_2+2z_3, z_3, z_1+z_1+2z_4)))_j\\
&\ot (-\max(-2z_2+z_3-z_4, -z_1-z_4, -z_2))_j
\ot(-\max(-2z_2+z_3, -z_1, z_3-2z_4))_i.
\end{align*}
We call these $\phi_{ij}^{(k)}$ \ ($k=0,1,2,3$) the {\it braid-type isomorphisms} of $B_i$'s. We omit $\phi_{ij}^{(3)}$ and $\phi_{ji}^{(3)}$, since we do not use in this article. Note that $\phi_{ij}^{(k)}$ and $\phi_{ji}^{(k)}$ are inverse to each other.
\end{pro}

By braid-type isomorphisms, for any $w\in W$ and its reduced words $i_1\cd i_l$ and 
$j_1\cd j_l$, we obtain
\[
B_{i_1}\ot\cd\ot B_{i_l}\cong B_{j_1}\ot\cd\ot B_{j_l}.
\]

\begin{df}\label{def of braid-moves}
We define a {\it braid-move} on the set of reduced words of $w\in W$
to be a composition of the following transformations 
induced from the braid relations:
\begin{eqnarray*}
&&\cd ij\cd\to \cd ji\cd (a_{ij}a_{ji}=0),\q
\cd iji\cd\to \cd jij\cd (a_{ij}a_{ji}=1),\\
&&\cd ijij\cd\to \cd jiji\cd (a_{ij}a_{ji}=2),\q
\cd ijijij\cd\to \cd jijiji\cd (a_{ij}a_{ji}=3),
\end{eqnarray*}
which are called 2-move, 3-move, 4-move, 6-move respectively.
\end{df}

\begin{df}[{\it Cellular crystal}]\label{def of cellular}
For a reduced word $\bfi=i_1i_2\cd i_k$ of $w \in W$, 
we call the crystal $\bbB_\bfi:=B_{i_1}\ot\cd\ot B_{i_k}$ 
a cellular crystal
associated with $\bfi$. 
\end{df}

For  a sequence of indices $\bfi=(i_1,\cd,i_k)\in I^k$, we identify $\bbB_{\bfi}$ with $\bbZ^k$ as follows:
\[
(x_1,\cd,x_m) \longmapsto \til f_{i_1}^{x_1}(0)_{i_1}
\ot\cd\ot \til f_{i_m}^{x_m}(0)_{i_m}=(-x_1)_{i_1}\ot\cd\ot(-x_m)_{i_m},
\]
where if $n<0$, then $\fit^n(0)_i$ means $\eit^{-n}(0)_i$.

By the tensor structure of crystals, we can describe 
the explicit crystal structure on $\bbB_\bfi:=B_{i_1}\ot \cd\ot B_{i_m}$ 
as follows: 
For $x=(x_1,\cd,x_m)\in \bbB_\bfi$, define
\begin{equation*}
 \sigma_k(x):=x_k+\sum_{j<k}\lan h_{i_k},\al_{i_j}\ran x_j
\label{sigma-def}
\end{equation*}
and for $i\in I$ define
\begin{eqnarray*}
&&\TY(\wtil\sigma,i,)(x):=\max\{\sigma_k(x)\,|\,1\leq k\leq m\,{\rm and}\,
i_k=i\},\label{simgak}\\
&&\TY(\wtil M,i,)=\TY(\wtil M,i,)(x):=\{k\,|\,1\leq k\leq m,\,i_k=i,\,
\sigma_k(x)=\TY(\wtil\sigma,i,)(x)\},
\label{Mi}
\\
&&
\TY(\wtil m,i,f)=\TY(\wtil m,i,f)(x):=\max\,\TY(\wtil M,i,)(x),\q
\TY(\wtil m,i,e)=\TY(\wtil m,i,e)(x):=\min\,\TY(\wtil M,i,)(x).
\label{mfme}
\end{eqnarray*}
The actions of the Kashiwara operators $\eit,\fit$ and the functions
$\vep_i,\vp_i$ 
and $\wt$  are written explicitly:
\begin{eqnarray}
&& \fit(x)_k:=x_k+\del_{k,\TY(\wtil m,i,f)},\qq\qq
\eit(x)_k:=x_k-\del_{k,\TY(\wtil m,i,e)},\\
&&
\wt(x):=-\sum_{k=1}^m x_k\al_{i_k},\q
\vep_i(x):=\TY(\wtil\sigma,i,)(x),\q
\vp_i(x):=\lan h_i,\wt(x)\ran+\vep_i(x).
\label{wt-vep-vp-ify}
\end{eqnarray}

\subsection{Procedure by braid moves}\label{procedure}
Let $w_0 \in W$ be the longest element. We fix  $\bfi=i_{l_1}i_{l_2}\cdots i_{i_l}$ as a reduced word of $w_0$. We present a procedure to move a given letter $i$ to the rightmost position of the fixed reduced word by successive braid moves. This procedure will be required in Section 5.

\subsubsection{Type $A_n$}
We adopt the following Dynkin diagram for type $A_n$:
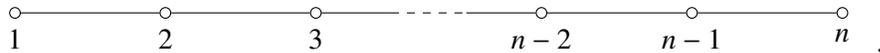
\[
\begin{tikzpicture}
\node[draw, circle, inner sep=0.5mm, label=below:{$1$}] (1) at (0,0) {};
\node[draw, circle, inner sep=0.5mm, label=below:{$2$}] (2) at (2,0) {};
\node[draw, circle, inner sep=0.5mm, label=below:{$3$}] (3) at (4,0) {};
\node[draw, circle, inner sep=0.5mm, label=below:{$n-2$}] (4) at (7,0) {};
\node[draw, circle, inner sep=0.5mm, label=below:{$n-1$}] (5) at (9,0) {};
\node[draw, circle, inner sep=0.5mm, label=below:{$n$}] (6) at (11,0) {};
\draw (node cs:name=1) --(node cs:name=2);
\draw (node cs:name=2) --(node cs:name=3);
\draw (3) -- (5,0);
\draw[dashed] (5,0) -- ++(1,0); 
\draw (6,0) -- (4);
\draw (node cs:name=4) --(node cs:name=5);
\draw (node cs:name=5) --(node cs:name=6); 
\node at (11.5,-0.5,0) {.};
\end{tikzpicture}
\]
Fix $\bfi=(1)(21) \cdots (n-1 \cdots 21)(n\cdots21)$ as a reduced longest word.
\begin{enumerate}
\item[{(Step1)}]
By applying 2-moves, move the $i$-th occurrence of the letter $1$ (reading from \textbf{right to left}) as far to the right as possible. Then we obtain the subword $121$ and apply a 3-move : $121 \mapsto 212$.
\item[{(Step2)}]
By applying 2-moves, move the right $2$ appeared in (Step1) as far to the right as possible. Then we obtain the subword $232$ and apply a 3-move : $232 \mapsto 323$.
\item[{(Step3)}]
Iterate (Step 2) successively for $3, 4, \dots, i-1$ until the subword $i(i-1)i$ appears.
\item[{(Step4)}]
Since there only exist letters strictly less than $i-1$ on the right side of the right $i$ appeared in (Step3), one can move the letter $i$ to the rightmost position by applying 2-moves.
\end{enumerate}

\begin{ex}\label{exA1}
Let $n=4$ and fix $ \bfi=1(21)(321)(4321) $ as a reduced longest word.
We move $3$ to the rightmost position as an example. The procedure in this case is
\begin{align*}
&12\textbf{1}3214321 \xrightarrow{(\text{Step1})} 123\textbf{1}214321 \xrightarrow{(\text{Step1})} 123\textbf{212}4321 \xrightarrow{(\text{Step2})} 123214\textbf{2}321\\
&\xrightarrow{(\text{Step2})} 123214\textbf{323}1 \xrightarrow{(\text{Step4})} 123214321\textbf{3} .
\end{align*}
\end{ex}

\subsubsection{Type$B_n$ and $C_n$}
We adopt the following Dynkin diagram for type $B_n$:
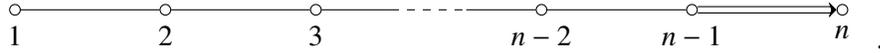
\[
\begin{tikzpicture}
\node[draw, circle, inner sep=0.5mm, label=below:{$1$}] (1) at (0,0) {};
\node[draw, circle, inner sep=0.5mm, label=below:{$2$}] (2) at (2,0) {};
\node[draw, circle, inner sep=0.5mm, label=below:{$3$}] (3) at (4,0) {};
\node[draw, circle, inner sep=0.5mm, label=below:{$n-2$}] (4) at (7,0) {};
\node[draw, circle, inner sep=0.5mm, label=below:{$n-1$}] (5) at (9,0) {};
\node[draw, circle, inner sep=0.5mm, label=below:{$n$}] (6) at (11,0) {};
\draw (1) --(2);
\draw (2) --(3);
\draw (3) -- (5,0);
\draw[dashed] (5,0) -- ++(1,0); 
\draw (6,0) -- (4);
\draw (4) --(5);
\draw[double distance=2pt, -{Stealth[length=3pt, width=6pt]}] (5)--(6);
\node at (11.5,-0.5,0) {.};
\end{tikzpicture}
\]
And for $C_n$, we adopt
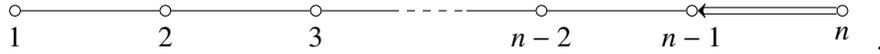
\[
\begin{tikzpicture}
\node[draw, circle, inner sep=0.5mm, label=below:{$1$}] (1) at (0,0) {};
\node[draw, circle, inner sep=0.5mm, label=below:{$2$}] (2) at (2,0) {};
\node[draw, circle, inner sep=0.5mm, label=below:{$3$}] (3) at (4,0) {};
\node[draw, circle, inner sep=0.5mm, label=below:{$n-2$}] (4) at (7,0) {};
\node[draw, circle, inner sep=0.5mm, label=below:{$n-1$}] (5) at (9,0) {};
\node[draw, circle, inner sep=0.5mm, label=below:{$n$}] (6) at (11,0) {};
\draw (1) --(2);
\draw (2) --(3);
\draw (3) -- (5,0);
\draw[dashed] (5,0) -- ++(1,0); 
\draw (6,0) -- (4);
\draw (4) --(5);
\draw[double distance=2pt, -{Stealth[length=3pt, width=6pt]}] (6)--(5);
\node at (11.5,-0.5,0) {.};
\end{tikzpicture}
\]
Fix $\bfi=(12 \cdots n)^n$ as a reduced longest word. We assume $1 \leq i \leq n-1$.
\begin{enumerate}
\item[{(Step1)}]
By applying 2-moves, move the $(i+1)$-th occurrence of 1 (reading from \textbf{left to right}) as far to the left as possible. Then we have the subword  $121$ and apply a 3-move : $121 \mapsto 212$.
\item[{(Step2)}]
By applying 2-moves, move the $2$ next to the 1 chosen in (Step1) as far to the left as possible. Then we have the subword $232$ and apply a 3-move : $232 \mapsto 323$.
\item[{(Step3)}]
Iterate (Step2) successively for $3,4, \dots, n-2$ until $(n-1)(n-2)(n-1)$ appears. Then we have $(n-1)n(n-1)n.$
\item[{(Step4)}]
Apply a 4-move : $(n-1)n(n-1)n \mapsto n(n-1)n(n-1)$.
\item[{(Step5)}]
By applying 2-moves, move the right $(n-1)$ appeared in (Step4) as far to the right as possible. Then we have $(n-1)(n-2)(n-1)$ and apply a 3-move : $(n-1)(n-2)(n-1) \mapsto (n-2)(n-1)(n-2)$.
\item[{(Step6)}]
By applying 2-moves, move the right $(n-2)$ appeared in (Step5) as far to the right as possible. Then we have $(n-2)(n-3)(n-2$ and make a 3-move : $(n-2)(n-3)(n-2) \mapsto (n-3)(n-2)(n-3)$.
\item[{(Step7)}]
Iterate (Step6) successively for $n-3, n-4, \dots, i+1$ until $i(i+1)i$ appears.
\item[{(Step8)}]
Since there only exist letters strictly greater than $i+1$ on the right side of the right $i$ appeared in (Step7), one can move the letter $i$ to the rightmost position by applying 2-moves.
\end{enumerate}
\medskip

\begin{ex}\label{eqB4}
Let n=4 and fix $ \bfi=(1234)(1234)(1234)(1234) $ as a reduced longest word. We move $2$ as an example. The procedure in this case is
\begin{align*}
&12341234\textbf{1}2341234 \xrightarrow{(\text{Step1})} 123412\textbf{1}342341234 \xrightarrow{(\text{Step1})} 1234\textbf{212}342341234\\
\xrightarrow{(\text{Step2})}&12342123\textbf{2}4341234 \xrightarrow{(\text{Step2})} 123421\textbf{323}4341234 \xrightarrow{(\text{Step4})}12342132\textbf{4343}1234\\
\xrightarrow{(\text{Step5})}&123421324341\textbf{3}234 \xrightarrow{(\text{Step5})}123421324341\textbf{232}4
\xrightarrow{(\text{Step8})}123421324341234\textbf{2}.
\end{align*}
\end{ex}

\subsubsection{Type $D_n$}
We adopt the following Dynkin diagram for type $D_n$:
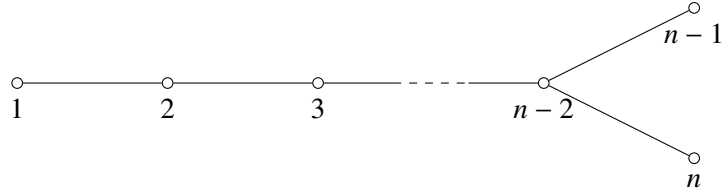
\[
\begin{tikzpicture}
\node[draw, circle, inner sep=0.5mm, label=below:{$1$}] (1) at (0,0) {};
\node[draw, circle, inner sep=0.5mm, label=below:{$2$}] (2) at (2,0) {};
\node[draw, circle, inner sep=0.5mm, label=below:{$3$}] (3) at (4,0) {};
\node[draw, circle, inner sep=0.5mm, label=below:{$n-2$}] (4) at (7,0) {};
\node[draw, circle, inner sep=0.5mm, label=below:{$n-1$}] (5) at (9,1) {};
\node[draw, circle, inner sep=0.5mm, label=below:{$n$}] (6) at (9,-1) {};
\draw (node cs:name=1) --(node cs:name=2);
\draw (node cs:name=2) --(node cs:name=3);
\draw (3) -- (5,0);
\draw[dashed] (5,0) -- ++(1,0); 
\draw (6,0) -- (4);
\draw (node cs:name=4) --(node cs:name=5);
\draw (node cs:name=4) --(node cs:name=6); 
\node at (9.5,-1.5,0) {.};
\end{tikzpicture}
\]
Fix $(1\cdots n)^{(n-1)}$ as a reduced longest word. We assume $1 \leq i \leq n-2$. 
\begin{enumerate}
\item[{(Step1)}]
By applying 2-moves, move the  $(i+1)$-th occurrence of the letter $1$ (reading from \textbf{left to right}) as far to the left as possible. Then we obtain the subword $121$ and apply a 3-move  $121\mapsto 212$.
\item[{(Step2)}]
By applying 2-moves, move the letter $2$ next to the letter $1$ chosen in (Step1) as far to the left as possible. Then we obtain the subword $232$ and apply a 3-move $232 \mapsto 323$.
\item[{(Step3)}]
Iterate (Step2) successively for $3, 4, \dots, n-3$ until $(n-2)(n-3)(n-2)$ appears. Then we obtain the subword $(n-2)(n-1)n(n-2)(n-1)n$.
\item[{(Step4)}]
Apply a 2-move $(n-2)(n-1)n(n-2)(n-1)n \mapsto (n-2)n(n-1)(n-2)(n-1)n$
\item[{(Step5)}]
Apply a 3-move: $(n-2)n(n-1)(n-2)(n-1)n \mapsto (n-2)n(n-2)(n-1)(n-2)n$. 
\item[{(Step6)}]
Apply a 3-move: $(n-2)n(n-2)(n-1)(n-2)n \mapsto n(n-2)n(n-1)(n-2)n$. 
\item[{(Step7)}]
Apply a 2-move: $n(n-2)n(n-1)(n-2)n \mapsto n(n-2)(n-1)n(n-2)n.$
\item[{(Step8)}]
Apply a 3 move: $n(n-2)(n-1)n(n-2)n \mapsto n(n-2)(n-1)(n-2)n(n-2).$
\item[{(Step9)}]
By applying 2-moves, move the right $(n-2)$ appeared in (Step9) as far to the right as possible. Then apply a 3 move : $(n-2)(n-3)(n-2) \mapsto (n-3)(n-2)(n-3).$
\item[{(Step10)}]
Iterate (Step9) successively for $n-3, n-4, \dots, i+1$ until the subword $i(i+1)i$ appears.
\item[{(Step11)}]
Since there only exist letters strictly greater than $i+1$ on the right side of the right $i$ appeared in (Step10), 
one can move the letter $i$ to the rightmost position by applying 2-moves.
\end{enumerate}

\begin{ex}\label{eq:D4}
Let $n=4$ and fix $\bfi=(1234)(1234)(1234)$ as a reduced longest word. We move $2$ as an example. The procedure in this case is
\begin{align*}
&12341234\textbf{1}234 \xrightarrow{(\text{Step1})} 123412\textbf{1}34234 \xrightarrow{(\text{Step1})} 1234\textbf{212}34234\\
\xrightarrow{(\text{Step4})}&1234212\textbf{43}234 \xrightarrow{(\text{Step5})} 12342124\textbf{232}4 \xrightarrow{(\text{Step6})}123421\textbf{424}324\\
\xrightarrow{(\text{Step7})}&12342142\textbf{34}24 \xrightarrow{(\text{Step8})}123421423\textbf{242}.
\end{align*}
\end{ex}

\section{Quiver Hecke Algebra and its modules}\label{QHA}
We review the basics of the quiver Hecke algebras (cf.\cite{K-L,Loc,jams,Rou}).
\subsection{Definition of Quiver Hecke Algebra}\label{defQHA}
For a finite index set $I$ and a field $\bf k$, 
let $({\mathscr Q}_{i,j}(u,v))_{i,j\in I}\in{\bf k}[u,v]$ be polynomials satisfying:
\begin{enumerate}
\item ${\mathscr Q}_{i,j}(u,v)={\mathscr Q}_{j,i}(v,u)$ for any $i,j\in I$.
\item ${\mathscr Q}_{i,j}(u,v)$ is in the form:
\[
{\mathscr Q}_{i,j}(u,v)=\begin{cases}\displaystyle
\sum_{a(\al_i,\al_i)+b(\al_j,\al_j)=-2(\al_i,\al_j)}
t_{i,j;a,b}u^av^b&\hbox{ if }i\ne j,\\
0&\hbox{ if }i=j,
\end{cases}
\]
where $t_{i,j;-a_{ij},0}\in{\bf k}^\times$.
\end{enumerate}
For $\beta=\sum_im_i\al_i\in Q_+$ with $|\beta|:=\sum_im_i=m$, set
$I^\beta:=\{\nu=(\nu_1,\cd,\nu_m)\in I^m\mid \sum_{k=1}^m\al_{\nu_k}=\beta\}$.
\begin{df}
For $\beta\in Q_+$, the {\it quiver Hecke algebra $R(\beta)$} associated with a Cartan matrix  $A$ and 
polynomials ${\mathscr Q}_{i,j}(u,v)$ is the $\bfk$-algebra generated by
\[
\{e(\nu)|\nu\in I^\beta\},\q
\{x_k|1\leq k\leq n\},\q
\{\tau_i|1\leq i\leq n-1\}
\]
with the following relations:
\begin{eqnarray*}
&&e(\nu)e(\nu')=
\del_{\nu,\nu'}e(\nu),\q
\sum_{\nu\in I^\beta}e(\nu)=1, \q
e(\nu)x_k=x_ke(\nu), \q
x_kx_l=x_lx_k,\\
&&\tau_le(\nu)=e(s_l(\nu))\tau_l,\q
\tau_k\tau_l=\tau_l\tau_k\,\,\hbox{ if }|k-l|>1,\\
&& \tau_k^2e(\nu)={\mathscr Q}_{\nu_k,\nu_{k+1}}(x_k,x_{k+1})e(\nu),\\
&&(\tau_kx_l-x_{s_k(l)}\tau_k)e(\nu)=\begin{cases}
-e(\nu)&\hbox{ if }l=k, \,\,\nu_k=\nu_{k+1},\\
e(\nu)&\hbox{ if }l=k+1, \,\,\nu_k=\nu_{k+1},\\
0&\hbox{otherwise},
\end{cases}\\
&&
(\tau_{k+1}\tau_k\tau_{k+1}-\tau_k\tau_{k+1}\tau_k)e(\nu)=
\begin{cases}
\ovl{\mathscr Q}_{\nu_k,\nu_{k+1}}(x_k,x_{k+1},x_{k+2})e(\nu)
&\hbox{ if } \nu_k=\nu_{k+2},\\
0&\hbox{ otherwise},
\end{cases}
\end{eqnarray*}
where $\ovl{\mathscr Q}_{i,j}(u,v,w)=\frac{{\mathscr Q}_{i,j}(u,v)-{\mathscr Q}_{i,j}(w,v)}{u-w}
\in\bfk[u,v,w]$.
\end{df}
Note that there exists an anti-automorphism $\psi$ on $R(\beta)$ that preserves all generators.
\begin{pro}
The relations above are homogeneous if we define 
\[
{\rm deg}(e(\nu))=0,\q
{\rm deg}(x_ke(\nu))=(\al_{\nu_k},\al_{\nu_k}),\q
{\rm deg}(\tau_l e(\nu))=-(\al_{\nu_l},\al_{\nu_{l+1}}).
\]
Thus, $R(\beta)$ becomes a $\bbZ$-graded algebra.
\end{pro}

\begin{df}
For $\beta,\gamma\in Q_+$, set $e(\beta,\gamma)
=\sum_{\nu\in I^\beta,\nu'\in I^\gamma}e(\nu,\nu')$. 
We define an injective homomorphism
$\xi(\beta,\gamma):R(\beta)\ot R(\gamma)\to 
e(\beta,\gamma)R(\beta+\gamma)e(\beta,\gamma)$ by 
\begin{align*}
&\xi(\beta,\gamma)(e(\nu)\ot e(\nu'))=e(\nu,\nu'),\\ 
&\xi(\beta,\gamma)(x_k e(\beta)\ot1)=x_ke(\beta,\gamma),\q
\xi(\beta,\gamma)(1\ot x_ke(\gamma))=x_{k+|\beta|}e(\beta,\gamma),\\
&\xi(\beta,\gamma)(\tau_ke(\beta)\ot1)=\tau_ke(\beta,\gamma),\q
\xi(\beta,\gamma)(1\ot\tau_ke(\gamma))=\tau_{k+|\beta|}e(\beta,\gamma).
\end{align*}
\end{df}
\begin{df}
Let $M=\bigoplus_{k\in\bbZ}M_k$ be a $\bbZ$-graded $R(\beta)$-module. Define
a {\it grading shift functor} $q$ on  
the category of graded $R(\beta)$-modules $R(\beta)$-Mod by 
\[
qM:=\bigoplus_{k\in\bbZ}(qM)_k,\q\hbox{where }
(qM)_k=M_{k-1}.
\]
For $M,\,N\in R(\beta)$-Mod, let
${\rm Hom}_{R(\beta)}(M,N)$ be the space of degree preserving morphisms and define 
$\textsc{Hom}_{R(\beta)}(M,N):=\bigoplus_{k\in \bbZ}{\rm Hom}_{R(\beta)}(q^kM,N)$, 
which is a space of morphisms up to grading shift. We define
${\rm deg}(f)=k$ for $f\in {\rm Hom}_{R(\beta)}(q^kM,N)$. 
\end{df}

\begin{df}
For $M\in R(\beta)$-Mod and $N\in R(\gamma)$-Mod, define 
the {\it convolution product} $\circ$ by 
\[
M\circ N:=R(\beta+\gamma)e(\beta,\gamma)\ot_{R(\beta)\ot R(\gamma)}
(M\ot N)
\]
$\circ$ is an exact bifunctor, that is, $\circ$ is an exact functor with respect to both variables. We also define 
\[
 M\nabla N:={\rm hd}(M\circ N),\ M\Del N:={\rm soc}(M\circ N), 
\]
where ${\rm hd}(M\circ N)$ is the head of $M\circ N$ and ${\rm soc}(M\circ N)$ is the socle of $M\circ N$.  
\end{df}
Recall that the head of a module $M$ is the quotient by its radical and the socle of a module $M$ is the summation of all simple submodules.

\begin{df}
Let $\psi$ be the anti-automorphism of $R(\beta)$ preserving all generators. 
For $M\in R(\beta)$-Mod, define $M^*:=\textsc{Hom}_\bfk(M,\bfk)$ with the $R(\beta)$-
module structure given by 
\[
(r\cdot f)(u):=f(\psi(r)u) \q \text{for}\ r\in R(\beta),\ u\in M\ \text{and}\ f\in M^*,
\] which is called a {\it dual module} of $M$.
In particular, if $M\cong M^*$ we call $M$ is {\it self-dual}.
\end{df}

\subsection{Categorification of quantum coordinate ring $\aq$}

Let $R(\beta)$-{\rm gmod} be the full subcategory of $R(\beta)$-Mod whose 
objects are  
finite-dimensional graded $R(\beta)$-modules and set
$R$-gmod$=\bigoplus_{\beta\in Q_+}R(\beta)$-gmod.
Define the functors 
\begin{align*}
E_i:R&(\beta)\hbox{-gmod}\to R(\beta-\al_i)\hbox{-gmod} \\
&M\longmapsto e(\al_i,\beta-\al_i)M, \\
F_i:R&(\beta)\hbox{-gmod}\to R(\beta+\al_i)\hbox{-gmod }\\
&M\longmapsto F_i(M)=\langle i \rangle \circ M.
\end{align*}
where $ e(\al_i,\beta-\al_i):=
\sum_{\nu\in I^\beta, \nu_1=i}e(\nu)$ and 
$\langle i \rangle:=R(\al_i)/R(\al_i)x_1$ is a 1-dimensional simple $R(\al_i)$-module. We also define $E_i^*$ and $F_i^*$ in the opposite manner.  Note that the functors $E_i$, $E_i^*$, $F_i$ and $F_i^*$ are exact functors.

Since the grading shift functor $q$ is exact, $\bbZ[q,q^{-1}]$ module structure is induced by $q[M]:=[qM]$. Moreover, since the convolution product $\circ$ is an exact bifunctor, multiplication on $\cK(R\hbox{\rm -gmod})$ is induced by $[M][N]:=[M\circ N]$. One can verify this multiplication satisfies associativity; hence, $\cK(R\hbox{\rm -gmod})$ becomes $\bbZ[q,q^{-1}]$-algebra.

The following theorem is one of the fundamental properties of the quiver Hecke algebra, which states that the quiver Hecke algebra categorifies the quantum group.
\begin{thm}[\cite{K-L,Rou}]
There exists an isomorphism of $\bbQ(q)$-algebras
\[
\bbQ(q)\ot_{\bbZ[q,q^{-1}}\cK(R\hbox{\rm -gmod})\cong \aq.
\]
\end{thm}

\subsection{Categorification of the crystal $B(\ify)$}
\begin{lem}[\cite{K-L}]\label{simple}
For any simple $R(\beta)$-module $M$, ${\rm soc}(E_i M)$, ${\rm hd}(E_i M)$
 and ${\rm hd}(F_i M)$ are all simple modules. Here we also have that
${\rm soc}(E_i M)\cong{\rm hd}(E_i M)$ up to grading shifts.
\end{lem}
For $M\in R(\beta)$-gmod, define
\begin{eqnarray}
&&\wt(M)=-\beta,\
\vep_i(M)=\max\{n\in\bbZ\,|\,E_i^nM\ne0\},\
\vp_i(M)=\vep_i(M)+\lan h_i,\wt(M) \ran,\\
&&\Eit M:= q_i^{1-\vep_i(M)}{\rm soc}(E_i M)\cong q_i^{\vep_i(M)-1}{\rm hd}(E_iM),\qq
\Fit M:=q_i^{\vep_i(M)}{\rm hd}(F_iM).\label{eitfit}
\end{eqnarray}
Set 
$\bbB:=\{S\,|\,S\hbox{ is a self-dual simple module in }R\hbox{-gmod}\}$.
Then, it follows from Lemma \ref{simple}
that $\wtil E_i$ and $\wtil F_i$ are well-defined on $\bbB$.
\begin{thm}[\cite{L-V}]\label{LV}
The 6-tuple $(\bbB,\{\Eit\},\{\Fit\},\wt,\{\vep_i\},\{\vp_i\})_{i\in I}$ holds 
a crystal structure and there exists  the following 
isomorphism of crystals:
\[
\Psi:\bbB\q\mapright{\sim}\q B(\ify).
\]
\end{thm}

\section{Localization}
We very briefly review the theory of the localization of monoidal categories. See  \cite{KKK}, \cite{Loc}, \cite{Loc2} and \cite{Loc3} for more details.

\subsection{Localization of monoidal caterogies via graded braider}\label{localization}
We recall the subcategory $\scC_w$ and its localization $\wtil\scC_w$.  For $M \in R(\beta)\hbox{-gMod}$, we define 
\begin{align*}
W(M)&:=\{\gamma\in Q_{+} \cap (\beta-Q_{+})\ | \ e(\gamma, \beta- \gamma)M \ne 0\}, 
\end{align*}
For $w\in W$, we define the full monoidal subcategory of $R\hbox{-gmod}$ by
\[
\scC_w:=\{M\in R\hbox{-gmod} \ | \ W(M) \subset Q_+ \cap wQ_-\}.
\]
We note that for finite type $\ge$ and the longest element $w=w_0$  in $W$, $\scC_{w_0}=\rgmod$.
For a Weyl group element $w$, let $s_{i_1}\cd s_{i_l}$ 
be its reduced expression. For a dominant weight $\Lm\in P_+$, 
set 
\[
m_k:=\lan h_{i_k},s_{i_{k+1}}\cd s_{i_l}\Lm\ran.\qq (k=1,\cd,l)
\] 
We define the
{\it determinantial module} associated with $w$ and $\Lm$ by
\[
{\bf M}(w\Lm,\Lm):=\wtil F_{i_1}^{m_1}\cd
\wtil F_{i_l}^{m_l}{\bf 1},
\]
where ${\bf 1}$ is a trivial $R(0)$-module. We have ${\bf M}(w\Lm, \Lm) \in \scC_w.$


Let $\Lm$ be a $\bbZ$-lattice, $( \cT=\oplus_{\lm\in \Lm} \cT_\lm, \ot)$ be a $\bfk$-linear $\Lm$-graded monoidal category, $q$ be the grading shift functor on $\cT$ and ${\bf 1}\in\cT_0$ be the unit object of $\cT$.
\begin{df}[\cite{Loc}]\label{real-com}
A {\it graded braider} is a triple
$(C,R_C,\phi)$, where $C\in \cT$, $\bbZ$-linear map $\phi:\Lm\to\bbZ$ and a 
morphism
\[
R_C:C\ot X\to q^{\phi(\lm)}X\ot C
\q(X\in\cT_\lm),
\]
which is functorial in $X \in \cT$ such that satisfying the following commutative diagram:
\[
\xymatrix{
\linethickness{30pt}
C\ot X\ot Y\ar@{->}^{R_C(X)\ot Y}[r]
\ar@{->}_{R_C(X\ot Y)}[dr]& q^{\phi(\lm)}\ot X\ot C\ot Y
\ar@{->}^{X\ot R_C(Y)}[d]\\
&q^{\phi(\lm+\mu)}(X\ot Y)\ot C,
}\qq
\xymatrix{
\linethickness{30pt}
C\ot \bf 1\ar@{->}^{R_C(\bf1)}[r]
\ar@{->}_{\cong}[dr]& {\bf1} \ot C
\ar@{->}^{\cong}[d]\\
&C.
}
\]
\end{df}

The triple 
\begin{align}\label{fam of braiders}
\{M(w\Lambda_i, \Lambda_i, R_{M(w\Lambda_i, \Lambda_i}, \phi_{M(w\Lambda_i, \Lambda_i} \}_{i\in I}
\end{align}
becomes a non-degenerate real commuting family of braiders. See \cite{Loc} for the definition of a real commuting family of braiders.

We apply the localization procedure to $\rgmod$ and $\scC_w$.
Let $\wtil\scC_w:=\scC_w[M(w\Lambda_i, \Lambda_i)^{\circ -1}; i\in I]$ denote the localization of $\scC_w$ by \eqref{fam of braiders}. We also let $\wtil\rgmod[w]:=\rgmod[M(w\Lambda_i, \Lambda_i)^{\circ -1}; i\in I]$ denote the localization of $\rgmod$ by the same graded braider. 
We set $C_\Lambda:=M(w\Lambda, \Lambda)$, in particuler,  $C_i:=M(w\Lambda_i, \Lambda_i)$. We summarize the propertise of $\wtil\rgmod[w]$.

\begin{pro}[\cite{Loc}]\label{til-pro}
Let $Q_w:\rgmod\to\tilRgmod[w]$ be the canonical functor. Then,
\begin{enumerate}
\item
$\tilRgmod[w]$ is an abelian category and the functor $Q_w$ is exact.
\item
For any simple object $S\in\rgmod$, $\ Q_w(S)$ is simple or $0$ in $\tilRgmod[w]$.
\item
$\wtil C_i:=Q_w(C_i)$ ($i\in I$) 
is invertible central graded braider in $\tilRgmod[w]$.
\medskip

\nd
For $\mu\in P$, define $\wtil C_\mu$ such that  
$\wtil C_\mu:=Q_w(C_\mu)$ for $\mu\in P_+$, 
$\wtil C_{-\Lm_i}=C_i^{\circ -1}$ and 
$\wtil C_{\lm+\mu}=\wtil C_{\lm}\circ \wtil C_{\mu}$ for $\lm,\mu\in P$
up to grading shifts.
\item
Any simple object in $\tilRgmod[w]$ is isomorphic to $\wtil C_\Lm\circ Q_w(S)$ 
for some simple module $S\in\rgmod$ and $\Lm\in P$.\label{5.12 (4)}
\end{enumerate}
Note that in (4) $\Lm\in P$ and $S\in\rgmod$ are not necessarily unique.
\end{pro}

For any $w \in W$, we have a quasi-commutative diagram
\[
\xymatrix@C=10em{
\scC_w \ \ar@{>->}[r]_{\iota_w} \ar[d]_{\Phi_w}& \rgmod \ar[d]_{Q_w} \\
\wtil\scC_w \ \ar@{>-->}[r]_{\widetilde\iota_w}& \wtil\rgmod [w].
}
\]
where $\Phi_w$ and $Q_w$ denote the localization functors  and $\widetilde \iota_w$ is a functor induced from  the inclusion functor $\iota_w$, which gives an equivalence of categories \cite{Loc3}. 
\subsection{Crystal Structure on localized quantum coordinate rings}
Following \cite{KN}, we define a crystal structure of $\wtil\rgmod[w]$, and introduce theorems that play a fundamental role in the next section. We first recall the notion of rigidity and the invariant $\Lambda$, which are required to define the crystal structure.

\begin{df}\label{rigid}
Let $X,Y$ be objects in a monoidal category $\cT$, and $\vep:X\ot Y\to 1$
and $\eta:1\to Y\ot X$ morphisms in $\cT$.
We say that a pair $(X,Y)$ is {\it dual pair} or 
$X$ is a {\it left dual} to  $Y$ or  
$Y$ is a {\it right dual} to  $X$ if 
the following compositions are identities:
\[
X\simeq X\ot 1\,\,\mapright{{\rm id}\ot \eta}\,\,X\ot Y\ot X
\,\,\mapright{\vep\ot {\rm id}}\,\,1\ot X\simeq X,\,\,
Y\simeq 1\ot Y\,\,\mapright{\eta\ot {\rm id}}\,\,Y\ot X\ot Y
\,\,\mapright{{\rm id}\ot \vep}\,\,Y\ot 1\simeq Y
\]
\end{df}
We denote a right dual to $X$ by $\scD(X)$ and  a 
left dual to $X$ by $\scD^{-1}(X)$.
\vspace{-5pt}
\begin{thm}[\cite{Loc}, \cite{Loc2}]
$\wtil\scC_w$ is {\it rigid}, i.e.,
every object in $\wtil\scC_w$ has left and right duals.
\end{thm}

\begin{df}[\cite{KN}]
Let $(\scC, \ot)$ be a graded monoidal category and $M,N \in \scC$ be simple objects.
When $\dim\textsc{Hom}(M\otimes N, N\otimes M)=1$, the pair $(M,N)$ is said to be {\it $\Lambda$-definable}. A non-zero morphism $\mathbf{r}_{M,N}$ in this space is called the {\it $R$-matrix}, and we set $\Lambda(M,N):=\mathrm{deg}(\mathbf{r}_{M,N}).$
\end{df}

Let $\mathrm{Irr}(\wtil\scC_w)$ be the set of equivalence classes of simple objects in $\wtil\scC_w$ up to grading shifts. 
We define a crystal structure on $\mathrm{Irr}(\wtil\scC_w)$.
\begin{thm}\cite{KN}
For $X \in \mathrm{Irr}(\wtil\scC_w)$, we define the operators by
\begin{align*}
\vep_i(X)&=d_i^{-1}\wtil\Lambda(Q_w(\lan i \ran) ,X),\ \vep_i^*=d_i^{-1}\wtil\Lambda(X, Q_w(\lan i \ran)),\\
\vp_i(X)&=\vep_i(X)+\lan h_i, \wt(X) \ran,\ \vp_i^*(X)=\vep_i^*(X)+\lan h_i, \wt(X) \ran,\\
\wtil F_i X&=q_i^{\vep_i(X)}Q_w(\lan i \ran) \nabla X,\ \wtil F_i^*=q_i^{\vep_i^*(X)}X \nabla Q_w(\lan i \ran),\\ 
\wtil E_i X&=q_i^{\vp_i(X)+1}X\nabla\scD Q_w(\lan i \ran),\ \wtil E_i^* X=q_i^{\vp_i^*(X)+1}\scD^{-1} Q_w(\lan i \ran)\nabla X,
\end{align*}
where $d_i=(\alpha_i, \alpha_i)/2$ and $\wtil\Lambda(M,N):=(\Lambda(M,N)+(\mathrm{wt}(M),\mathrm{wt}(N)))/2$. Then these operators define a crystal structure on $\mathrm{Irr}(\wtil\scC_w)$.
\end{thm}

Let $\bfi=i_1 i_2\cdots i_l$ be a reduced word, $w=s_{i_1}s_{i_2}\cdots s_{i_l}$ be a reduced expression of $w$ and $\bbB_\bfi$ be the cellular crystal associated to $\bfi$. For $M \in \scC_w$,  we define $K_{\bfi}(M) \in \bbZ^l$ by 
\begin{align*}
&M_l=M,\ c_k=\vep_{i_k}^*(M_k) \ \text{for} \ 1 \leq k \leq l,\\
&M_{k-1}=(\wtil E_{i_k}^*)^{c_k}(M_k),\ K_{\bfi}(M)=(c_1,\dots, c_l).
\end{align*}
We regard $K_{\bfi}(M)$ as an element of $\bbB_\bfi$ by
$(c_1, \dots, c_l) \mapsto \wtil f_{i_1}^{c_1} (0)_{i_1}\ot\cdots\ot\wtil f_{i_l}^{c_l} (0)_{i_l}$. 
We extend $K_{\bfi} : \scC_w \rightarrow \bbZ^l$ to $K_{\bfi} : \wtil \scC_w \rightarrow \bbZ^l$ by
\[
\wtil C_{\Lambda}^{-1}\circ Q_w(M) \mapsto K_{\bfi}(M)-K_{\bfi}(C_{\Lambda})
\] 
for $\Lambda \in P_+$ \cite{KN}.

The following theorem is necessary to define $\vep_i^*$ on $\bbB_\bfi$,
\begin{thm}[\cite{KN}]\label{crisom2}
Let $w \in W$ and $i \in I$ satisfying $ w'=ws_i < w$. Define $\bfE_{w',w}^* : \mathrm{Irr}(\wtil\scC_w)\longrightarrow \mathrm{Irr}(\wtil\scC_{w'})$ by 
\[
(\wtil C_{\Lambda})^{-1}\circ Q_w(M) \mapsto Q_{w'}(\wtil E_i^{* \mathrm{max}}C_{\Lambda})^{-1}\circ Q_{w'}( \wtil E_i^{* \mathrm{max}}M).
\]
Then we have a bijective map
\[
\Psi_{w',w}:\mathrm{Irr}(\wtil\scC_w)\longrightarrow \mathrm{Irr}(\wtil\scC_{w'})\times \bbZ \quad (X \mapsto (\bfE_{w',w}^*(X), \vep_i^*(X))
\]
and an isomorphism of crystals $K_{\bfi} : \mathrm{Irr}(\wtil\scC_w) \rightarrow \bbZ^l \simeq \bbB_\bfi$.
\end{thm}

%
\section{Characterization of the unit object}
As an application of Theorem \ref{crisom2}, we give a characterization of the unit object of $\wtil\rgmod[w]$. In the rest of this article, we restrict ourselves to the case where $\ge$ is of classical finite type. Let $w=w_0$ is the longest element of the associated Weyl group $W$ and $\bfi$ be a reduced word associated to $w_0$, which we call a {\it  reduced longest word}. Since the fixed reduced word $\bfi$ is the longest one, any letter $i$ can be moved to the rightmost position of the reduced longest word by applying the braid moves properly, and each braid move induces the braid type isomorphism $\phi_{ij}^{(k)}$ \ ($k=0,1,2,3$). Thus, this procedure induces an isomorphism of the crystals 
\[
\phi : \bbB_\bfi \to \bbB_{\mathbf{i'}}\ot B_i \quad (x \mapsto \phi(x)=b'\ot (a)_i)
\]
where $\phi$ is the composition of the braid type isomorphisms $\phi_{ij}^{(k)}$, and $\bfi$ and $\mathbf{i'}i$ are the reduced words of $w$. For $x \in \bbB_\bfi$, we then set $\vep_{i}^*(x)=-a$. We may describe $\vep_{i}^*(x)$ to be the integer indexed by $i$ when the letter $i$ is moved to the rightmost position of the reduced longest word. Note that this definition of  $\vep_{i}^*(x)$ is motivated by $\Psi_{w',w}$ in Theorem \ref{crisom2}, and $\vep_i^*(x)$ depends neither on the choice of the reduced longest word nor on a process of braid-moves (cf. \cite{Kana-N}). We sometimes write $\vep_{i}^*$ for $\vep_{i}^*(x)$ when there is no risk of confusion.

 We express  $x \in \bbB_\bfi$ as
\[
x=(z_{1,i_{1}})_{i_{1}}\ot\dots\ot(z_{j, i_{k}})_{i_{k}}\ot\dots\ot(z_{s,i_{n}})_{i_{n}},
\]
where the double index $(j, i_{k})$ of $z_{j, i_K}$ indicates the $j$-th occurrence of $i_{k}$ reading from left to right. For example, for $\bfi=1(21)(321)$, which is a longest word of $A_3$, we express
\[
x=(z_{1,1})_1\ot(z_{1,2})_2\ot(z_{2,1})_1\ot(z_{1,3})_3\ot(z_{2,2})_2\ot(z_{3,1})_1.
\]
\subsection{Main results}
We now provide the explicit form of $\vep_i^*(x)$ for a fixed longest word $\bfi$ for each type and state the main result of this article.
\begin{pro}\label{pro-vep*}
Let W be the Weyl group, $\bfi$ = $i_1\cd i_l$ \ be a reduced longest word and $\bbB_\bfi$ be the cellular crystal associated $\bf{i}$. In what follows, we understand $z_{j,k}=0$ whenever $k=0$.  
\begin{enumerate}
\item
Type $A_n$ :  Fix $\bfi=(1)(21) \dots (n-1 \dots 21)(n\dots21)$ as a reduced longest word. For $x \in \bbB_\bfi$ expressed as
\[
x=(z_{1,1})_1\ot(z_{1,2})_2\ot(z_{2,1})_1\ot\cdots\ot(z_{1,n})_n\ot\cdots\ot(z_{n-1,2})_2\ot(z_{n,1})_1,
\]
$\vep_{i}^*(x)$ is given by
\begin{align*}
&\vep_{i}^*(x)=\underset{1 \leq k \leq i}{\max}(z_{n-i+2, k-1}- z_{n-i+1, k}) \q (1\leq i \leq n).
\end{align*}
\item
Type $B_n$ : Fix $\bfi=(12 \dots n)^n$ as a reduced longest word. For
$x \in \bbB_\bfi$ expressed as
\[
x=(z_{1,1})_1\ot(z_{1,2})_2\ot\cdots(z_{1,n})_n\ot(z_{2,1})_1\ot\cdots\ot(z_{n,n-1})_{n-1}\ot(z_{n,n})_n,
\]
$\vep_{i}^*(x)$ is given by
\begin{align*}
&\vep_{i}^*=\underset{1 \leq k \leq n-i-1}\max(-\zeta_i, z_{i+k+1,n-k}-z_{i+k+1,n-k-1}) \q (1 \leq i \leq n-2), \\ 
&\zeta_i=-\max(-z_{i+1,n-1}+z_{i+1,n}, -\eta_i, z_{i+1,n-1}-z_{i,n}) \q (1 \leq i \leq n-1),\\
&\eta_i=-\underset{1 \leq k \leq n-1}\max(z_{i+1,k-1}-z_{i,k}) \q (i \leq i \leq n-1),\\
&\vep_{n-1}^*=-\zeta_{n-1}, \ \vep_{n}^*=-z_{n.n}. 
\end{align*}
\item
Type $C_n$ : Fix $\bfi=(12 \dots n)^n$ as a reduced longest word. $x \in \bbB_\bfi$ has the same expression as in type $B_n$.  $\vep_{i}^*(x)$ and $\eta_i(x)$ are the same as those in type $B_n$, while $\zeta_i$ is given by
\begin{align*}
\zeta_i(x)=-\max(-2z_{i,n}+z_{i+1,n-1}, -\eta_i, z_{i+1,n-1}-2z_{i+1,n}) \ (1 \leq i \leq n-1).
\end{align*}
\item
Type $D_n$ : Fix $\bfi=(12 \dots n)^{n-1}$ as a reduced longest word. For $x \in \bbB_\bfi$ expressed as
\[
x=(z_{1,1})_1\ot(z_{1,2})_2\ot\cdots(z_{1,n})_n\ot(z_{2,1})_1\ot\cdots\ot(z_{n-1,n-1})_{n-1}\ot(z_{n-1,n})_n,
\]
$\vep_{i}^*(x)$ is given by 
\begin{align*}
\vep_{i}^* 
&= \underset{1 \leq k \leq n-i-2}\max(-\kappa_i, z_{i+k+1,n-k-1} - z_{i+k+1,n-k-2}) \q (1 \leq i \leq n-3), \\
\kappa_i 
&= -\max\Bigg( 
\begin{array}{lcr}
-\theta_i,\ z_{i+1,n-1} - z_{i,n},\ z_{i+1,n-2} - z_{i,n-1} - z_{i,n}, \\
 z_{i+1,n} - z_{i,n-1},\ z_{i+1,n} + z_{i+1,n-1} - z_{i+1,n-2}
\end{array}
\Bigg) \q  (1 \leq i \leq n-2), \\
\theta_i 
&= -\underset{1 \leq k \leq n-2}{\max}(z_{i+1, k-1} - z_{i,k})
\q (1 \leq i \leq n-2), \\
&\vep_{n-2}^*=-\kappa_{n-2},\ \vep_{n-1}^*=-z_{n-1, n-1},\ \vep_{n}^*=-z_{n-1,n}. 
\end{align*}
\end{enumerate}
\end{pro}
Using these forms of $\vep_i^*$, we obtain the following result.
\begin{thm}\label{mainthm}
Let W be the Weyl group of classical finite type.  \ $\bfi$ = $i_1\cd i_l$ \ be a reduced longest word, and $\bbB_\bfi$ be the cellular crystal associated $\textbf{i}$. For $x \in \bbB_\bfi$, the following are equivalent.
\begin{enumerate}
\item
$\mathrm{wt}(x)=0$ and for all $i \in I, \vep_{i}^*(x)=0.$
\item
$x= (0)_{i_1}\otimes \cd \otimes(0)_{i_l}\in \bbB_\bfi.$
\end{enumerate}
\end{thm}

We will prove Proposition \ref{pro-vep*} and Theorem \ref{mainthm} in the subsequent subsections, treating each type separately. Before that, we state an immediate consequence of Theorem \ref{mainthm}, which provides a characterization of the unit object in $\wtil\rgmod[w]$.
\begin{cor}
For $X\in \mathrm{Irr}(\wtil\scC_w)$, the following are equivalent.
\begin{enumerate}
\item
$\wt(X)=0$ and for all $i \in I, \vep_{i}^*(X)=0.$
\item
$X=\mathbf1.$
\end{enumerate}
\end{cor}
{\sl Proof.}
This follows from the bijectivity of $K_{\bfi}$.
\qed

\subsection{Proof of Proposition \ref{pro-vep*} and Theorem \ref{mainthm} for Type $A_n$}
First, we will prove Proposition \ref{pro-vep*} for type $A_n$. As explained above, the procedure  in Section \ref{procedure} for type $A_n$. induces an isomorphism 
\[
\phi(x)=x'\ot(-\vep_i^*(x))_i. \q (x \in \bbB_{\bfi}, x' \in \bbB_{\bfi'})
\] 
We show that the explicit form of ${\vep_{i}^*}(x)$ induced by the procedure is given by
\begin{align}\label{eqAn}
\vep_{i}^*(x)=\underset{1 \leq k \leq i}{\max}(z_{n-i+2, k-1}- z_{n-i+1, k}) \q  (1 \le i \le n) .
\end{align}
We proceed by induction on $n$, the rank of $\ge$. For $n=1$, the claim is obvious. Assume that \eqref{eqAn} is true for $n-1$. Let $\bfi'= (1)(21)\cdots(n-1\cdots1)$ be a longest word of $A_{n-1}$, then $\bfi=\bfi'(n\cdots1)$. Applying the procedure for type $A_n$ in Section \ref{procedure} to the subword $\bfi'$ of $\bfi$, we have $\bfi=\bfi'(n\cdots1) \mapsto \bfi''(i-1)(n\cdots1)$, that is, the letter $i-1$ in $\bfi'$ is moved to the rightmost position. This induces a crystal morphism $\bbB_\bfi \to \bbB_{\bfi''}\ot B_{i-1}\ot B_n\ot B_{n-1}\ot\cdots\ot B_1$,
\[
x \mapsto b''\ot(-\vep_{i-1}^{*(n-1)})_{i-1}\ot(z_{1,n})_n\ot(z_{2,n-1})_{n-1}\cdots(z_{n,1})_1
\]
where $b''\in \bbB_{\bfi''}$. By the induction hypothesis, 
\begin{align*}
\vep_{i-1}^{* \ (n-1)}&=\underset{1\leq k \leq i-1}\max(z_{(n-1)-(i-1)+2,k-1}-z_{(n-1)-(i-1)+1,k})\\
&=\underset{1\leq k \leq i-1}\max(z_{n-i+2,k-1}-z_{n-i+1,k}).
\end{align*}
Applying  2-moves, we have the subword $(i-1)i(i-1)$. Since $(i-1)$ appears $n-i+1$ times  and $i$ appears $n-i$ times in $\bfi'$, the corresponding part of the element takes the form
\[
\cdots\ot(-\vep_{i-1}^{* \ (n-1)})_{i-1}\ot(z_{n-i+1,i})_i\ot(z_{n-i+2,i-1})_{i-1}\ot\cdots.
\]
Applying a 3-move, we have the subword $i(i-1)i$, which induces the braid type isomorphism $\phi_{i-1i}^{(1)}$, and the corresponding part of the element becomes
\[
\cdots\ot(A)_i\ot(B)_{i-1}\ot(-\max(\vep_{i-1}^{* \ (n-1)}, z_{n-i+2,i-1}-z_{n-i+1,i})_i \cdots
\]
for some integers A and B. Therefore, we conclude
\begin{align*}
\vep_i^{*} &=\max(-\vep_{i-1}^{* \ (n-1)}, z_{n-i+2,i-1}-z_{n-i+1,i})\\
&=\max(\underset{1\leq k \leq i-1}\max(z_{n-i+2,k-1}-z_{n-i+1,k}), z_{n-i+2,i-1}-z_{n-i+1,i})\\
&=\underset{1 \leq k \leq i}{\max}(z_{n-i+2, k-1}- z_{n-i+1, k}).
\end{align*}
Thus the claim is true for $n$ as well.
\qed

\begin{ex}
Recall Example \ref{exA1}. The algorithm induces the following composition of morphisms:
\begin{align*}
&(z_{1,1})_1\ot(z_{1,2})_2\ot(z_{2,1})_1\ot(z_{1,3})_3\ot(z_{2,2})_2\ot(z_{3,1})_1\ot(z_{1,4})_4\ot(z_{2,3})_3\ot(z_{3,2})_2\ot(z_{4,1})_1 \\ 
\mapsto&(z_{1,1})_1\ot(z_{1,2})_2\ot(z_{1,3})_3\ot (z_{2,1})_1\ot(z_{2,2})_2\ot(z_{3,1}) _1\ot(z_{1,4})_4\ot(z_{2,3})_3\ot(z_{3,2})_2\ot(z_{4,1})_1\\
\mapsto&(z_{1,3})_3\ot(X)_2\ot(Y)_1\ot(-\max(-z_{2,1},z_{3,1}-z_{2,2}))_2\ot(z_{1,4})_4\ot(z_{2,3})_3\ot(z_{3,2})_2\ot(z_{4,1})_1\\
\mapsto&(z_{1,3})_3\ot(X)_2\ot(Y)_1\ot(z_{1,4})_4\ot(-\max(-z_{2,1},z_{3,1}-z_{2,2}))_2\ot(z_{2,3})_3\ot(z_{3,2})_2\ot(z_{4,1})_1\\
\mapsto&(S)_3\ot(T)_2\ot(-\max(\max(-z_{2,1},z_{3,1}-z_{2,2}), z_{3,2}-z_{2.3}))_3\ot(z_{4,1})_1\\
\mapsto&(S)_3\ot(T)_2\ot(z_{4,1})_1\ot(-\max(\max(-z_{2,1},z_{3,1}-z_{2,2}),z_{3,2}-z_{2.3})_3.\\
\end{align*}
 Hence, we have
\[\vep_{3}^*=\max(\max(-z_{2,1},z_{3,1}-z_{2,2}),z_{3,2}-z_{2.3})=\max(-z_{2,1},z_{3,1}-z_{2,2},z_{3,2}-z_{2.3}).\]
\end{ex}

We now prove Theorem \ref{mainthm} for type $A_n$.
By Proposition \ref{pro-vep*}, we have
\begin{align}
&\vep_{i}^*=\underset{1 \leq k \leq i} \max(z_{n-i+2,k-1}-z_{n-i+1,k})=0 \implies z_{n-i+2,k-1} \leq z_{n-i+1,k} \label{A_n-1} \\ &(1 \leq k \leq i ,\ 2 \leq i \leq n), \nn\\
&\vep_{1}^*=-z_{n,1}=0, \label{A_n-2}\\
&\mathrm{wt}(x)=0 \implies \sum_{j=1}^{n-i+1}z_{j,i}=0 \q (1 \leq i \leq n).\label{A_n-3}
\end{align}
 By \eqref{A_n-1} with $k=i$, we have \ $z_{n-i+2, i-1}, \leq z_{n-i+1,i}$ \ $(2 \leq i \leq n)$, which together with \eqref{A_n-2} and \eqref{A_n-3}, implies
\[
0=z_{n,1} \leq z_{n-1,2} \leq \cdots \leq z_{2,n-1} \leq z_{1,n}=0.
\]
Thus, we obtain $z_{i,j}=0$ for all $i, j$ such that $i+j=n+1$. Since\ $z_{1,n-1}+z_{2,n-1}=0$ by \eqref{A_n-3}, it follows that \ $z_{1,n-1}=0.$

Next, by \eqref{A_n-1} with $k=i-1$, we have \ $z_{n-i+2, i-2} \leq z_{n-i+1,i-1}$ \ $(3 \leq i \leq n)$. Hence, together with \ $z_{1,n-1}=0$, we obtain 
\[
0 \leq z_{n-1,1} \leq \cdots \leq z_{2, n-2} \leq z_{1, n-1}=0,
\]
which implies $z_{i,j}=0$ for all $i,j$ such that $i+j=n$. Since $z_{1,n-2}+z_{2,n-2}+z_{3,n-2}=0$ and  $z_{2,n-2}=z_{3,n-2}=0$, we have $z_{1,n-2}=0$.
Repeating the same argument, we obtain $z_{i,j}=0$ for all $i,j$.
\qed
\subsection{Proof of Proposition \ref{pro-vep*} and Theorem \ref{mainthm} for Type $B_n$ and $C_n$}
First, we prove Proposition \ref{pro-vep*} for types $B_n$ and $C_n$, that is,  we show that the function $\vep_i^*$ induced by the procedure described in Section \ref{procedure} for types $B_n$ and $C_n$ is given by
\begin{align*}
&\vep_{i}^*=\underset{1 \leq k \leq n-i-1}\max(-\zeta_i, z_{i+k+1,n-k}-z_{i+k+1,n-k-1}) \q (1 \leq i \leq n-2),  \\
&\zeta_i=-\max(-z_{i+1,n-1}+z_{i+1,n}, -\eta_i, z_{i+1,n-1}-z_{i,n}) \q \text{for type $B_n$},\\
&\zeta_i=-\max(-2z_{i,n}+z_{i+1,n-1}, -\eta_i, z_{i+1,n-1}-2z_{i+1,n}) \q \text{for type $C_n$} \\
& (1 \leq i \leq n-1), \nn \\
&\eta_i=-\underset{1 \leq k \leq n-1}\max(z_{i+1,k-1}-z_{i,k}) \q (i \leq i \leq n-1),\\
&\vep_{n-1}^*=-\zeta_{n-1}, \ \vep_{n}^*=-z_{n.n}.
\end{align*}
We employ the following notation:
\begin{align*}
&\eta_{i,1}:=z_{i,1},\ \eta_{i,j}:=-\max(-\eta_{i,j-1}, z_{i+1,j-1}-z_{i,j})=-\underset{1 \leq k \leq j}\max( z_{i+1,k-1}-z_{i,k}) \q ( 2 \leq j \leq n-1),\\
&\vep_{i,1}^*:=-\max(-\zeta_i, z_{i+2,n-1}-z_{i+2,n-2}), \\
&\vep_{i,j}^*:=-\max(-\vep_{i,j-1}, z_{i+j+1,n-j}-z_{i+j+1,n-j-1})=-\underset{1 \leq k \leq j}\max(-\zeta_i, z_{i+k+1,n-k}-z_{i+k+1,n-k-1})\\
& (2 \leq j \leq n-i-1).
\end{align*}
By definition, we have $\eta_{i}=\eta_{i,n-1}, \vep_{i}^*=-\vep_{i, n-i-1}^*$.
 (Step1) induces the following isomorphisms
\begin{align*}
&\cdots\ot(z_{i,1})_1\ot(z_{i,2})_2\ot\cdots\ot(z_{i+1,1})_1\ot\cdots \\
\mapsto&\cdots\ot(z_{i,1})_1\ot(z_{i,2})_2\ot(z_{i+1,1})_1\ot\cdots \\
\mapsto&\cdots\ot(A)_2\ot(B)_1\ot(-\max(-z_{i,1}, z_{i+1,1}-z_{i,2}))_2\ot\cdots \\
=&\cdots\ot(A)_2\ot(B)_1\ot(\eta_{i, 2})_2\ot\cdots
\end{align*}
where A and B are some integers. In the following, capital letters denote some integers. (Step2) yields
\begin{align*}
&\cdots\ot(\eta_{i, 2})_2\ot(z_{i,3})_3\ot\cdots\ot(z_{i+1, 2})_2\ot\cdots \\
\mapsto&\cdots\ot(\eta_{i, 2})_2\ot(z_{i,3})_3\ot(z_{i+1, 2})_2\ot\cdots \\
\mapsto&\cdots\ot(C)_3\ot(D)_2\ot(-\max(-\eta_{i,2}, z_{i+1,2}-z_{i,3}))_3\ot\cdots \\
=&\cdots\ot(C)_3\ot(D)_2\ot(\eta_{i,3})_3\ot\cdots
\end{align*}
in the corresponding part of the element. 
Thus, repeating this recursive procedure until $(n-1)(n-2)(n-1)$ appears, which is (Step3), yields
\begin{align*}
&\cdots\ot(\eta_{i,n-2})_{n-2}\ot(z_{i,n-1})_{n-1}\ot(z_{i,n})_n\ot(z_{i+1,n-2})_{n-2}\ot\cdots \\
\mapsto&\cdots\ot(\eta_{i,n-2})_{n-2}\ot(z_{i,n-1})_{n-1}\ot(z_{i+1,n-2})_{n-2}\ot(z_{i,n})_n\ot\cdots \\
\mapsto&\cdots\ot(E)_{n-1}\ot(F)_{n-2}\ot(-\max(\eta_{i,n-2}, z_{i+1, n-2}-z_{i,n-1}))_{n-1}\ot(z_{i,n})_n\ot\cdots \\
=&\cdots\ot(E)_{n-1}\ot(F)_{n-2}\ot(\eta_i)_{n-1}\ot(z_{i,n})_n\ot\cdots
\end{align*}
Therefore, the procedure (Step1) through (Step3) yields $\eta_i$. Now we must have the subword $(n-1)n(n-1)n$ and 
\[
\cdots\ot(\eta_i)_{n-1}\ot(z_{i,n})_n\ot(z_{i+1,n-1})_{n-1}\ot(z_{i+1, n})_{n}\ot\cdots
\]
in the corresponding part of the element. Then (Step4) induces
\begin{align*}
&\cdots\ot(\eta_i)_{n-1}\ot(z_{i,n})_n\ot(z_{i+1,n-1})_{n-1}\ot(z_{i+1, n})_{n}\ot\cdots \\
\mapsto&\cdots\ot(G)_{n}\ot(H)_{n-1}\ot(I)_{n}\ot(-\max(-z_{i+1,n-1}+z_{i+1,n}, -\eta_i, z_{i+1,n-1}-z_{i,n}))_{n-1}\ot\cdots\\
=&\cdots\ot(G)_{n}\ot(H)_{n-1}\ot(I)_{n}\ot(\zeta_i)_{n-1}\ot\cdots.
\end{align*}
Hence, $\zeta_i$ corresponds to (Step4). If $i=n-1$, we have $\vep_{n-1}^*=\zeta_{n-1}$ and so we assume $ i < n-1$ hereafter. We show the remaining part of the procedure yields $\vep_{i}^*$. (Step5) induces the following composition of isomorphisms
\begin{align*}
&\cdots\ot(\zeta_i)_{n-1}\ot\cdots\ot(z_{i+2,n-2})_{n-2}\ot(z_{i+2,n-1})_{n-1}\ot\cdots \\
\mapsto&\cdots\ot(\zeta_i)_{n-1}\ot(z_{i+2,n-2})_{n-2}\ot(z_{i+2,n-1})_{n-1}\ot\cdots \\
\mapsto&\cdots\ot(J)_{n-2}\ot(K)_{n-1}\ot(-\max(-\zeta_i, z_{i+2,n-1}-z_{i+2,n-2}))_{n-2}\ot\cdots \\
=&\cdots\ot(J)_{n-2}\ot(K)_{n-1}\ot(\vep_{i, 1}^*)_{n-2}\ot\cdots.
\end{align*}
Similarly, (Step6) induces 
\begin{align*}
&\cdots\ot(\vep_{i, 1}^*)_{n-2}\ot\cdots\ot(z_{i+3, n-3})\ot(z_{i+3, n-2})_{n-2}\ot\dots\\
\mapsto&\cdots\ot(\vep_{i, 1}^*)_{n-2}\ot(z_{i+3, n-3})\ot(z_{i+3, n-2})_{n-2}\ot\dots\\
\mapsto&\cdots\ot(L)_{n-3}\ot(M)_{n-2}\ot(-\max(\vep_{i, 1}^* ,z_{i+3, n-2}- z_{i+3, n-3}))_{n-3}\ot\cdots \\
=&\cdots\ot(L)_{n-3}\ot(M)_{n-2}\ot(\vep_{i,2}^*)_{n-3}\ot\cdots.
\end{align*}
From (Step7) we have 
\begin{align*}
&\cdots\ot(\vep_{i,n-i-2})_{i+1}\ot\cdots\ot(z_{n,i})_i\ot(z_{n, i+1})_{i+1}\ot\cdots \\
\mapsto&\cdots\ot(\vep_{i,n-i-2})_{i+1}\ot(z_{n,i})_i\ot(z_{n, i+1})_{i+1}\ot\cdots \\
\mapsto&\cdots\ot(N)_{i}\ot(O)_{i+1}\ot(-\max(-\vep_{i,n-i-2}, z_{n,i+1}-z_{n,i}))_{i}\ot\cdots\\
=&\cdots\ot(N)_{i}\ot(O)_{i+1}\ot(\vep_{i,n-i-1}^*)_{i}\ot\cdots.
\end{align*}
Thus, we conclude
\begin{align*}
&\vep_{i}^*=-\vep_{i,n-i-1}^*=\underset{1 \leq k \leq n-i-1}\max(-\zeta_i, z_{i+k+1,n-k}-z_{i+k+1,n-k-1}) \q (1 \leq i \leq n-2), \\
&\vep_{n-1}^*=\zeta_{n-1}.
\end{align*}
This proves Proposition \ref{pro-vep*} for type $B_n$. For type $C_n$, since $a_{n,n-1}=-1, a_{n-1,n}=-2$, applying $\phi_{n-1,n}^{(2)}$ in (Step4) gives the desired expression of $\zeta_{i}$ for type $C_n$.
\qed

\begin{ex}
Recall Example \ref{eqB4}. Then the algorithm induces the following composition of morphisms:
\begin{align*}
&(z_{2,1})_1\ot(z_{2,2})_2\ot(z_{2,3})_3\ot(z_{2,4})_4\ot(z_{3,1})_1\ot(z_{3,2})_2\ot(z_{3,3})_3\ot(z_{3,4})_4\ot\cdots\\
\mapsto&(z_{2,1})_1\ot(z_{2,2})_2\ot(z_{3,1})_1\ot(z_{2,3})_3\ot(z_{2,4})_4\ot(z_{3,2})_2\ot(z_{3,3})_3\ot(z_{3,4})_4\ot\cdots\\
\mapsto&\cdots\ot(-\max(-z_{2,1}, z_{3,1}-z_{2,2}))_2\ot(z_{2,3})_3\ot(z_{2,4})_4\ot(z_{3,2})_2\ot(z_{3,3})_3\ot(z_{3,4})_4\ot\cdots\\
\mapsto&\cdots\ot(-\max(-z_{2,1}, z_{3,1}-z_{2,2}))_2\ot(z_{2,3})_3\ot(z_{3,2})_2\ot(z_{2,4})_4\ot(z_{3,3})_3\ot(z_{3,4})_4\ot\cdots\\
\mapsto&\cdots\ot(-\max(\max(-z_{2,1}, z_{3,1}-z_{2,2}), z_{3,2}-z_{2,3})_3\ot(z_{2,4})_4\ot(z_{3,3})_3\ot(z_{3,4})_4\ot\cdots\\
\mapsto&\cdots\ot(A)_4\ot(B)_3\ot(C)_4\ot(-\max(-z_{3,3}+z_{3,4}, \eta_{2}, z_{3,3}-z_{2,4}))_3\ot(z_{4,1})_1\ot(z_{4,2})_2\ot(z_{4,3})_3\ot(z_{4,4})_4\\
\mapsto&\cdots\ot(z_{4,1})_1\ot(\zeta_{2})_3\ot(z_{4,2})_2\ot(z_{4,3})_3\ot(z_{4,4})_4\\
\mapsto&\cdots\ot(D)_2\ot(E)_3\ot(-\max(-\zeta_2, z_{4,3}-z_{4,2}))_2\ot(z_{4,4})_4\\
\mapsto&\cdots\ot(D)_2\ot(E)_3\ot(z_{4,4})_4\ot(-\max(-\zeta_2, z_{4,3}-z_{4,2}))_2.\\
\end{align*}
Hence, we have
\[
\vep_{2}^{*}=\max(-\zeta_2, z_{4,3}-z_{4,2}),\ \zeta_2=-\max(z_{3,3}+z_{3,4}, \eta_{2}, z_{3,3}-z_{2,4}),  \eta_2=-\max(-z_{2,1}, z_{3,1}-z_{2,2}).
\]
\end{ex}

{\sl Proof of Theorem \ref{mainthm} for type $B_n$ and $C_n$.}
We first prove Theorem \ref{mainthm} for type $B_n$. By Proposition \ref{pro-vep*}, we have
\begin{align}
&\vep_{i}^*=0 \implies -\zeta_i \leq 0, \ z_{i+k+1,n-k}-z_{i+k+1,n-k-1} \leq 0 \q (1 \leq k \leq n-i-1), \label{eq:bn1}\\
&-\zeta_i \leq 0 \implies  -\eta_i \leq 0,\label{eq:bn2} \\
&-\eta_i \leq 0 \implies z_{i+1,k-1}-z_{i,k} \leq 0, \ -z_{i,1}\leq 0 \q (1 \leq k \leq n-1) \label{eq:bn3}\\
&(1 \leq i \leq n-1),\nn\\
&\mathrm{wt}(x)=0 \implies \sum_{j=1}^{n}z_{j,i}=0 \q (1 \leq k \leq n),\label{eq:bn4}\\
&\vep_n^*=z_{n,n}=0.\label{eq:bn5}
\end{align}
From \eqref{eq:bn1} with $k=n-i-1$ and \eqref{eq:bn5}, we have 
\begin{equation}
0=z_{n,n}\leq \cdots \leq z_{n,i+1} \leq z_{n, i} \leq \cdots \leq z_{n,1}.
\end{equation}
By \eqref{eq:bn3} we have $z_{i,1} \geq 0 \ (1 \leq i \leq n-1)$. \eqref{eq:bn4} implies $z_{n,1} \leq 0$. Hence, we obtain $\ z_{n,i}=0\ (1 \leq i \leq n).$ 
Since $z_{i,1} \geq 0$ and \ $\sum_{j=1}^{n-1}z_{j,1}=0,$ it must be $\ z_{i,1}=0 \ (1 \leq i \leq n).$
By \eqref{eq:bn3} and \ $z_{i,1}=0$, we have 
\begin{equation}
0 = z_{i,1} \leq z_{i-1,2}\leq \dots \leq z_{1,i} \q (1 \leq i \leq n). \label{eq:bn6}
\end{equation}
\eqref{eq:bn6} implies $z_{i,j}\geq0$ for $i, j$ such that  $i+j\leq n+1.$
By \eqref{eq:bn3} and \ $z_{n,i}=0$, we have
\begin{equation}
0 = z_{n,i} \leq z_{n-1,i+1} \leq \dots \leq z_{i,n} \q (1 \leq i \leq n). \label{eq:bn7}
\end{equation}
\eqref{eq:bn7} implies $z_{i,j}\geq0$ for $i, j$ such that  $i+j\geq n+1$.
Combining, we obtain $\ z_{i,j} \geq 0$ for all $i, j$. It follows from \eqref{eq:bn4} that $\ z_{i,j} = 0$.

For type $C_n$, the only difference lies in the expression of $\zeta_i$, and the difference does not affect the proof. Hence, the same argument applies to type $C_n$.
\qed
\subsection{Proof of Propositon \ref{pro-vep*} and Theorem \ref{mainthm} for Type $D_n$}
We show Proposition \ref{pro-vep*} for type $D_n$, that is, we show $\vep_i^*$ induced by the procedure in Section \ref{procedure} for type $D_n$ is given by
\begin{align*}
\vep_{i}^* 
&= \underset{1 \leq k \leq n-i-2}\max(-\kappa_i, z_{i+1+k,n-k-1} - z_{i+1+k,n-k-2})\q (1 \leq i \leq n-3), \\
\kappa_i 
&= -\max\Bigg( 
\begin{array}{lcr}
-\theta_i,\ z_{i+1,n-1} - z_{i,n},\ z_{i+1,n-2} - z_{i,n-1} - z_{i,n}, \\
 z_{i+1,n} - z_{i,n-1},\ z_{i+1,n} + z_{i+1,n-1} - z_{i+1,n-2}
\end{array}
\Bigg), \\
\theta_i 
&= -\underset{1 \leq k \leq n-2}{\max}(z_{i+1, k-1} - z_{i,k})
\q  (1 \leq i \leq n-2),\\
&\vep_{n-2}^*=-\kappa_{n-2},\ \vep_{n-1}^*=-z_{n-1, n-1},\ \vep_{n}^*=-z_{n-1,n}. 
\end{align*}
We set
\begin{align*}
&\theta_{i,1}:=-z_{i,1},\ \theta_{i,j}:=-\max(-\theta_{i,j-1}, z_{i, j-1}-z_{i,j}) -\underset{1 \leq k \leq j}{\max}(z_{i+1, j-1} - z_{i,j})\q (2 \leq j \leq n-2), \\
&\vep_{i,1}^*=-\max(-\kappa_i, z_{i+1, n-2}-z_{i+1,n-3}), \\
&\vep_{i,j}^*=-\max(-\vep_{i,j-1}^*, z_{i+1, n-1}-z_{i+1,n-2})=-\underset{1 \leq k \leq n-i-2}\max(-\kappa_i,\ {\max}(z_{i+1+k,n-k-1} - z_{i+1+k,n-k-2})) \\
&(2 \leq j \leq n-i-2).
\end{align*}
By definition, we have $\theta_{i}=\theta_{i,n-2}$ and $\vep_{i}^*=-\vep_{i,n-i-2}$. Since (Step1), (Step2) and (Step3) of the procedure for type $D_n$ are almost the same as those for type $B_n$, those steps yield
\begin{align*}
&\cdots\ot(\theta_{i,n-3})_{n-3}\ot(z_{i,n-2})_{n-2}\ot(z_{i,n-1})_{n-1}\ot(z_{i,n})_n\ot(z_{i+1,n-3})_{n-3}\ot\cdots \\
\mapsto&\cdots\ot(\theta_{i,n-3})_{n-3}\ot(z_{i,n-2})_{n-2}\ot(z_{i+1,n-3})_{n-3}\ot(z_{i,n-1})_{n-1}\ot(z_{i,n})_n\ot\cdots \\
\mapsto&\cdots\ot(A)_{n-2}\ot(B)_{n-3}\ot(-\max(-\theta_{i,n-3}, z_{i+1, n-3}-z_{i+1, n-2})_{n-2}\ot(z_{i,n-1}))_{n-1}\ot(z_{i,n})_n\ot\cdots \\
=&\cdots\ot(A)_{n-2}\ot(B)_{n-3}\ot(\theta_{i})_{n-2}\ot(z_{i,n-1}))_{n-1}\ot(z_{i,n})_n\ot\cdots. 
\end{align*}
Now the part of the element takes the form of
\[
\cdots\ot(\theta_{i})_{n-2}\ot(z_{i,n-1})_{n-1}\ot(z_{i,n})_n\ot(z_{i+1,n-2})_{n-2}\ot(z_{i+1,n-1})_{n-1}\ot(z_{i+1,n})_n\ot\cdots.
\]
Recall that the procedure from (Step4) to (Step8) consists of the following successive braid moves:
\begin{align*}
&(n-2)(n-1)n(n-2)(n-1)n\xrightarrow{(\text{Step4})}(n-2)n(n-1)(n-2)(n-1)n\\
\xrightarrow{\text{(Step5)}}&(n-2)n(n-2)(n-1)(n-2)n\xrightarrow{\text{(Step6)}}n(n-2)n(n-1)(n-2)n\\
\xrightarrow{\text{(Step7)}}&n(n-2)(n-1)n(n-2)n\xrightarrow{\text{(Step8)}}n(n-2)(n-1)(n-2)n(n-2).
\end{align*}
Then these braid moves induce the following composition of isomorphisms:
\begin{align*}
&\cdots\ot(\theta_i)_{n-2}\ot(z_{i,n-1})_{n-1}\ot(z_{i,n})_n\ot(z_{i+1,n-2})_{n-2}\ot(z_{i+1,n-1})_{n-1}\ot(z_{i+1,n})_n\ot\cdots\\
\xrightarrow{(\text{Step4})}&\cdots\ot(\theta_i)_{n-2}\ot(z_{i,n})_n\ot(z_{i,n-1})_{n-1}\ot(z_{i+1, n-2})_{n-2}\ot(z_{i+1,n-1})_{n-1}\ot(z_{i+1,n})_n\ot\cdots\\
\xrightarrow{\text{(Step5)}}&
\begin{array}{r@{\quad}l}
\begin{gathered}
\cdots\ot(\theta_i)_{n-2}\ot(z_{i,n})_n\ot(\max(z_{i+1,n-1}, z_{i+1,n-2}-z_{i,n-1}))_{n-2}\\\ot
(z_{i,n-1}+z_{i+1,n-1})_{n-1}\ot(-\max(-z_{i,n-1}, z_{i+1,n-1}-z_{i+1,n-2}))_{n-2}\ot(z_{i+1,n})_n\ot\cdots
\end{gathered}
\end{array}\\
\xrightarrow{\text{(Step6)}}&
\begin{array}{r@{\quad}l}
\begin{gathered}
\cdots\ot(C)_{n}\ot(D)_{n-2}\ot(-\max(-\theta_i, z_{i+1,n-1}-z_{i,n}, z_{i+1,n-2}-z_{i,n-1}-z_{i.n})_{n}\\\ot
(z_{i,n-1}+z_{i+1,n-1})_{n-1}\ot(-\max(-z_{i,n-1}, z_{i+1,n-1}-z_{i+1,n-2}))_{n-2}\ot(z_{i+1,n})_n\ot\cdots
\end{gathered}
\end{array}\\
\xrightarrow{\text{(Step7)}}&
\begin{array}{r@{\quad}l}
\begin{gathered}
\cdots\ot(z_{i,n-1}+z_{i+1,n-1})_{n-1}\ot
(-\max(-\theta_i, z_{i+1,n-1}-z_{i,n}, z_{i+1,n-2}-z_{i,n-1}-z_{i.n})_{n}\\
\ot(-\max(-z_{i,n-1}, z_{i+1,n-1}-z_{i+1,n-2}))_{n-2}\ot(z_{i+1,n})_n\ot\cdots
\end{gathered}
\end{array}\\
\xrightarrow{\text{(Step8)}}&\cdots\ot(C)_{n-2}\ot(D)_n\ot(\kappa_i)_{n-2}\ot\cdots.
\end{align*}
Hence, (Step4) through (Step8) correspond to $\kappa_i$. If $i=n-2$, we obtain $\vep_{n-2}^*=\kappa_{n-2}$, therefore we assume $i < n-2$ in what follows. We show that the remaining part of  the procedure yields $\vep_i^*$, which is similar to the case of type $B_n$.
From (Step9), we have
\begin{align*}
&\cdots\ot(\kappa_i)_{n-2}\ot\cdots\ot(z_{i+2, n-3})_{n-3}\ot(z_{i+2, n-2})_{n-2}\ot\cdots\\
\mapsto&\cdots\ot(\kappa_i)_{n-2}\ot(z_{i+2, n-3})_{n-3}\ot(z_{i+2, n-2})_{n-2}\ot\cdots\\
\mapsto&\cdots\ot(\ot(E)_{n-3}\ot(F)_{n-2}\ot(-\max(-\kappa_i, z_{i+2, n-2}-z_{i+2, n-3}))_{n-3}\ot\cdots\\
=&\cdots\ot(\ot(E)_{n-3}\ot(F)_{n-2}\ot(\vep_{i,1}^*)_{n-3}\ot\cdots.
\end{align*}
Similarly, (Step10) yields
\begin{align*}
&\cdots\ot(\vep_{i,1}^*)_{n-3}\ot\cdots\ot(z_{i+3, n-4})_{n-4}\ot(z_{i+3, n-3})_{n-3}\ot\cdots\\
\mapsto&\cdots\ot(\vep_{i,1}^*)_{n-3}\ot(z_{i+3, n-4})_{n-4}\ot(z_{i+3, n-3})_{n-3}\ot\cdots\\
\mapsto&\cdots\ot(G)_{n-3}\ot(H)_{n-2}\ot(-\max(-\vep_{i,1}^*, z_{i+3, n-3}-z_{i+3, n-4}))_{n-4}\ot\cdots\\
=&\cdots\ot(G)_{n-3}\ot(H)_{n-2}\ot(\vep_{i,2}^*)_{n-4}\ot\cdots.
\end{align*}
Inductively, from (Step11) we obtain
\begin{align*}
&\cdots\ot(\vep_{i,n-i-3}^*)_{i+1}\ot\cdots\ot(z_{n-1, i})_{i}\ot(z_{n-1, i+1})_{i+1}\ot\cdots\\
\mapsto&\cdots\ot(\vep_{i,n-i-3}^*)_{i+1}\ot(z_{n-1, i})_{i}\ot(z_{n-1, i+1})_{i+1}\ot\cdots\\
\mapsto&\cdots\ot(J)_{i}\ot(K)_{i+1}\ot(-\max(-\vep_{i,n-i-3}^*, z_{n-1, i+1}-z_{n-1, i}))_{i}\ot\cdots\\
=&\cdots\ot(G)_{i}\ot(H)_{i+1}\ot(\vep_{i,n-i-2}^*)_{i}\ot\cdots.
\end{align*}
Thus, we conclude
\begin{align*}
&\vep_i^*=-\vep_{i,n-i-2}^*= \underset{1 \leq k \leq n-i-2}\max(-\kappa_i, z_{i+k+1,n-k-1} - z_{i+k+1,n-k-2}) \q  (1 \leq i \leq n-3),\\
&\vep_{n-2}^*=-\kappa_{n-2}.
\end{align*}
This proves Proposition \ref{pro-vep*} for type $D_n$.
\qed

\begin{ex}
Recall Example \ref{eq:D4}. The procedure induces the following composition of morphisms:
\begin{align*}
&\cdots(z_{2,1})_1\ot(z_{2,2})_2\ot(z_{2,3})\ot(z_{2,4})_4\ot(z_{3,1})_1\ot(z_{3,2})_2\ot(z_{3,3})_3\ot(z_{3,4})_4\\
\mapsto&\cdots(z_{2,1})_1\ot(z_{2,2})_2\ot(z_{3,1})_1\ot(z_{2,3})_3\ot(z_{2,4})_4\ot(z_{3,2})_2\ot(z_{3,3})_3\ot(z_{3,4})_4\\
\mapsto&\cdots\ot(-\max(-z_{2,1}, z_{3,1}-z_{2,2}))_2\ot(z_{2,3})_3\ot(z_{2,4})_4\ot(z_{3,2})_2\ot(z_{3,3})_3\ot(z_{3,4})_4\\
=&\cdots\ot(\theta_2)_2\ot(z_{2,3})_3\ot(z_{2,4})_4\ot(z_{3,2})_2\ot(z_{3,3})_3\ot(z_{3,4})_4\\
\mapsto&\cdots\ot(\theta_2)_2\ot(z_{2,4})_4\ot(z_{2,3})_3\ot(z_{3,2})_2\ot(z_{3,3})_3\ot(z_{3,4})_4\\
\mapsto&\cdots\ot(\theta_2)_2\ot(z_{2,4})_4\ot(\max(z_{3,3}, z_{3,2}-z_{2,3}))_2\ot(C)_3\ot(-\max(-z_{2,3}, z_{3,3}-z_{3,2}))_2\ot(z_{3,4})_4\\
\mapsto&\cdots\ot(-\max(-\theta_2, z_{3,3}-z_{2,4}, z_{3,2}-z_{2,3}-z_{2,4}))_4\ot(C)_3\ot(-\max(-z_{2,3}, z_{3,3}-z_{3,2}))_2\ot(z_{3,4})_4\\
\mapsto&\cdots\ot(C)_3\ot(-\max(-\theta_2, z_{3,3}-z_{2,4}, z_{3,2}-z_{2,3}-z_{2,4}))_4\ot(-\max(-z_{2,3}, z_{3,3}-z_{3,2}))_2\ot(z_{3,4})_4\\
\mapsto&\cdots\ot(D)_2\ot(E)_4\ot(-\max(-\theta_2, z_{3,3}-z_{2,4}, z_{3,2}-z_{2,3}-z_{2,4}, z_{3,4}-z_{2,3}, z_{3,4}+z_{3,3}-z_{3,2}))_2,\\
=&\cdots\ot(D)_2\ot(E)_4\ot(\kappa_2)_2.
\end{align*}
Hence, we have
\begin{align*}
&\vep_{2}^{*}=-\kappa_2=\max(-\theta_2, z_{3,3}-z_{2,4}, z_{3,2}-z_{2,3}-z_{2,4}, z_{3,4}-z_{2,3}, z_{3,4}+z_{3,3}-z_{3,2}),\\
&\theta_2=-\max(-z_{2,1}, z_{3,1}-z_{2,2}).
\end{align*}
\end{ex}

{\sl Proof of Theorem \ref{mainthm} for type $D_n$.}
 By Proposition \ref{pro-vep*}, we have
\begin{align}
&\vep_i^*=0 \implies z_{i+k+1, n-k-1} \leq z_{i+k+1, n-k-2}, 
\ -\kappa_i \leq 0 \q (1\leq k \leq n-i-2)\label{eq:dn1} \\
&(1 \leq i \leq n-3), \nn \\
&-\kappa_i \leq 0 \implies z_{i+1,n-1} \leq z_{i,n},\ z_{i+1,n} \leq z_{i,n-1},\ z_{i+1,n} + z_{i+1,n-1} \leq z_{i+1,n-2}, -\theta_i \leq 0, \label{eq:dn2}\\
&-\theta_i \leq 0 \implies z_{i+1, j-1} \leq z_{i,j} \q (1 \leq  j \leq n-2) \label{eq:dn3}\\
&(1 \leq i \leq n-2), \nn \\
&\mathrm{wt}(x)=0 \implies \sum_{j=1}^{n-1}z_{j,i}=0 \q (1 \leq i \leq n), \label{eq:dn4}\\
&\vep_{n-1}^*=z_{n-1, n-1}=0,\ \vep_{n}^*=z_{n-1,n}=0. 
\end{align}
Since $z_{n-1, n-1}, z_{n-1,n}=0$, by \eqref{eq:dn2} with $i=n-2$, we obtain
\begin{align*}
0=z_{n-1,n}+z_{n-1, n-1} \leq z_{n-1,n-2}.
\end{align*}
Recall the convention\ $z_{i, 0}=0$. By \eqref{eq:dn3} with $j=1$, we have $ 0 \leq z_{i,1} \ (1\leq i \leq n-2)$, and $z_{n-1,1}\leq 0$ by \eqref{eq:dn4}. If $z_{n-1,1}=0$, it follows that  $z_{i,1}=0 \ (1 \leq i \leq n-1)$ from  \eqref{eq:dn4}.
By \eqref{eq:dn1} with $k=n-i-2$, we have $z_{n-1, i+1}\leq z_{n-1, i} \ (1 \leq i \leq n-3)$, that is,
\[
0=z_{n-1, n-2}\leq \cdots \leq z_{n-1,1} \leq 0.
\]
Thus, we obtain $z_{n-1,1}=0$ and $z_{i,1}=0 \ (1 \leq i \leq n-1)$. By \eqref{eq:dn3} with $j=2$, we have $0=z_{i+1, 1} \leq z_{i,2} \ (1 \leq i \leq n-2)$. Since $z_{n-1, 2}=0$, we obtain $z_{i, 2}=0 \ (1 \leq i \leq n-1)$ by \eqref{eq:dn4}. Therefore, by induction, we have 
\[
z_{i,j}=0 \q (1 \leq i \leq n-1,\ 1\leq j \leq n-2).
\]
Let us show $z_{i, n-1}=z_{i, n}=0 \ (1 \leq i \leq n-1).$ By \eqref{eq:dn2} with $i=l, l-1$, we have 
\[
z_{l+1, n-1} \leq z_{l,n} \leq z_{l-1, n-1},\ z_{l+1,n} \leq z_{l, n-1} \leq z_{l-1, n} \q (2 \leq l \leq n-2).
\]
Since $0=z_{n-1, n-1}=z_{n-1,n}$, these imply $z_{l,n-1}$ and  $z_{l, n}$ are all non-negative for any $l$. Therefore, by \eqref{eq:dn4} we conclude $z_{i, n-1}=z_{i, n}=0 \ (1 \leq i \leq n-1)$. \qed

\bibliographystyle{amsalpha}

\begin{thebibliography}{A}

%
%
%
%
%
           


\bibitem{Kac} Kac V.G., 
Infinite dimensional Lie algebras 3rd ed., Cambridge Univ.Press, (1990).

\bibitem{Kana-N}
 Kanakubo Y. and Nakashima T., Half potential on geometric 
 crystals and connectedness of cellular crystals, Transformation Groups {\bf 28}, (2023), no. 1, 327–373.

\bibitem{KKK} Kang S-J., Kashiwara M. and Kim M., 
Symmetric quiver Hecke algebras and R-matrices of quantum affine algebras, Invent. Math. {\bf 211}, 
(2018), no. 2, 591--685.

\bibitem{jams} Kang S-J.,Kashiwara M., Kim M. and Oh S-J.,
Monoidal categorification of cluster algebras, J.Amer.Math.Soc., 
{\bf 31}, (2017), no. 2, 349--426.

\bibitem{simple}\bysame,
Simplicity of heads and
socles of tensor products, Compos. Math. {\bf 151}, (2015), no. 2, 377--396.

\bibitem{K1} Kashiwara M.,
 On crystal bases of the $q$-analogue of universal enveloping algebras,
	Duke Math. J., {\bf 63}, (1991), no. 2, 465--516.


\bibitem{K3}\bysame,
      Crystal base and Littelmann's refined Demazure character formula.
	     Duke Math. J. {\bf 71}, (1993), no. 3, 839--858.

\bibitem{KN} Kashiwara M. and Nakashima T., Crystal structure of localized quantum unipotent coordinate category, arXiv:2502.14319.

\bibitem{KP} Kashiwara M. and Park E., Affinizations and R-matrices for quiver Hecke algebras, J. Eur. Math.Soc. {\bf 20}, (2018), no. 5, 1161--1193.

\bibitem{strat}
Kashiwara M, Kim M.,  Oh S.-j., and Park E., 
Monoidal categories associated with strata of flag
manifolds, Adv. Math. {\bf 328} (2018), 959--1009.

\bibitem{Loc}\bysame,
Localization  for quiver 
 Hecke algebras, Pure Appl.Math. Q. {\bf 17}, (2021), 
no. 4, 1465--1548.

\bibitem{Loc2}\bysame,
Localizations for quiver Hecke algebras II. Proc. Lond. Math. Soc. (3) {\bf 127}, (2023), no. 4, 1134--1184.


\bibitem{Loc3}\bysame,
Localizations for quiver Hecke algebras III. Math. Ann. {\bf 390}, (2024), no. 4, 5075–-5108. 


\bibitem{K-L} Khovanov M. and Lauda A., A diagrammatic approach to categorification of quantum groups I, Represent. Theory {\bf 13}, (2009), 309--347.

\bibitem{K-L2}\bysame,
A diagrammatic approach to categorification of quantum groups II, Trans. Amer. Math. Soc. {\bf 363}, (2011), 2685--2700.   


\bibitem{L-V}  Lauda A.D. and Vazirani M., Crystals from 
categorified quantum groups, Adv.Math., {\bf 228}, (2011), 803--861.
%
%
%
%
\bibitem{N1} Nakashima T., Polyhedral Realizations of 
         Crystal Bases and Braid-type Isomorphisms, 
           Contemporary Mathematics {\bf 248}, (1999), 419--435.
%
%
%

\bibitem{N}\bysame,
Categorified Crystal Structure on Localized Quantum Coordinate Rings, arXiv:2208.08396v2.
%

\bibitem{Rou} Rouquier R., 2-Kac-Moody algebras, arXiv:0812.5023v1.

\bibitem{Rou2}\bysame,  Quiver Hecke algebras and 2-Lie algebras, Algebra Colloq. {\bf19}, (2012), no. 2, 359--410.



\end{thebibliography}

\end{document}